\documentclass{amsart}
\usepackage[latin1]{inputenc}
\usepackage[all]{xy}
\usepackage{amssymb}
\usepackage{rotating}
\usepackage{amsmath}
\usepackage{graphicx}
\usepackage{epsfig}
\DeclareMathOperator{\rank}{\rm{rank}}
\DeclareMathOperator{\tors}{\rm{tors}}
\DeclareMathOperator{\discr}{\rm{discr}}
\DeclareMathOperator{\Tr}{\rm{Tr}}
\DeclareMathOperator{\MW}{\rm{MW}}
\DeclareMathOperator{\Aut}{\rm{Aut}}
\DeclareMathOperator{\Km}{\rm{Km}}
\newtheorem{theorem}{Theorem}[section]
\newtheorem{rem}{Remark}[section]
\newtheorem{prop}{Proposition}[section]
\newtheorem{lemma}{Lemma}[section]

\newtheorem{defi}{Definition}[section]

\newcommand{\bprf}{{\it Proof.~}}
\newcommand{\binf}{{\it In fact }}
\newcommand{\ra}{\rightarrow}
\newcommand{\eprf}{\hfill $\square$ \smallskip\par}
\newcommand{\erem}{\hfill $\square$}
\newcommand{\PP}{ \mathbb{P}}
\newcommand{\C }{ \mathbb{C}}

\newcommand{\Z}{\mathbb{Z}}
\newcommand{\Q}{\mathbb{Q}}
\newcommand{\cl}{\mathcal{L}}

\input epsf
\input cyracc.def

\makeatletter
\def\blfootnote{\xdef\@thefnmark{}\@footnotetext}

\begin{document}

%\date{\today}

%\input prepictex
%\input pictex
%\input postpictex
\title{Elliptic fibrations and symplectic automorphisms\\ on K3 surfaces}
\author{Alice Garbagnati and Alessandra Sarti}
\address{Alice Garbagnati, Dipartimento di Matematica, Universit\`a di Milano,
  via Saldini 50, I-20133 Milano, Italy}

\email{alice.garbagnati@unimi.it}

\address{Alessandra Sarti, Institut f\"ur  Mathematik, Universit\"at Mainz,
Staudingerweg 9, 55099 Mainz, Germany. {\it Current address}: Universit\'e de Poitiers,
Laboratoire de Math\'ematiques et Applications,
 T\'el\'eport 2
Boulevard Marie et Pierre Curie
 BP 30179,
86962 Futuroscope Chasseneuil Cedex, France}

\email{sarti@math.univ-poitiers.fr}
\urladdr{http://www-math.sp2mi.univ-poitiers.fr/~sarti/}

\begin{abstract}
Nikulin has classified all finite abelian groups acting symplectically on a K3 surface and he has shown that the induced action on the  K3 lattice $U^3\oplus E_8(-1)^2$ depends only on the group but not on the K3 surface. For all the groups in the list of Nikulin we compute the invariant sublattice and its orthogonal complement by using some special elliptic K3 surfaces.
 \end{abstract}

\maketitle

%\tableofcontents
\pagestyle{myheadings}
\markboth{Alice Garbagnati and Alessandra Sarti}{Elliptic fibrations and symplectic automorphisms on K3 surfaces }


\section{Introduction}
An automorphism of a K3 surface is called symplectic if the
induced action on the holomorphic 2-form is trivial. The finite
groups acting symplectically on a K3 surface are classified by
\cite{Nikulin symplectic}, \cite{mukai}, \cite{gang}, where also
their fixed locus is described. In \cite{Nikulin symplectic}
Nikulin shows that the action of a finite abelian group of
symplectic automorphisms on the K3 lattice $U^3\oplus E_8(-1)^2$
is unique (up to isometries of the lattice), i.e. it depends only
on the group and not on the K3 surface. Hence one can consider a
special K3 surface, compute the action, then, up to isometry, this
is the same for any other K3 surface with the same finite abelian
group of symplectic automorphisms. It turns out that elliptic K3
surfaces are good candidates to compute this action, in fact one
can produce symplectic automorphisms by using sections.\\ The
finite abelian groups in the list of Nikulin are the following
fourteen groups:
\begin{eqnarray*}\begin{array}{l}\Z/n\Z,\ 2\leq n\leq 8,\
(\Z/m\Z)^2,\ m=2,3,4,\\ \Z/2\Z\times \Z/4\Z,\ \Z/2\Z\times
\Z/6\Z,\ (\Z/2\Z)^i,\ i=3,4.\end{array}\end{eqnarray*} For all the
groups $G$  except $(\Z/2\Z)^i,\ i=3,4$, there exists an elliptic
K3 surface with symplectic group of automorphisms $G$ generated by
sections of finite order (cf. \cite{shimada}). In general it is
difficult to describe explicitly  the action on the K3 lattice. An
important step toward this identification is to determine the
invariant sublattice, its orthogonal complement and the action of
$G$ on this orthogonal complement. In \cite{Nikulin symplectic}
Nikulin gives only rank and discriminant of these lattices,
however the discriminants are the correct ones only in the case that
the group is cyclic (cf. also \cite{alinikulin}). In the case of
$G=\Z/2\Z$ the action on the K3 lattice was computed by Morrison
(cf. \cite{morrison}) and in the cases of $G= \Z/p\Z,\,p=3,5,7$ we
computed in \cite{symplectic prime} the invariant sublattice and
its orthogonal complement in the K3 lattice. In this paper we
conclude the description of these lattices for all the fourteen
groups, in particular we can compute their discriminants, which are
not always the same as those given in \cite{Nikulin symplectic}.
It is interesting that some of the orthogonal complements to the
invariant lattices are very well known lattices. We denote by
$\Omega_G:=(H^2(X,\Z)^G)^{\perp}$, then

$$
\begin{array}{l|lllll}
G&\Z/3\Z &(\Z/2\Z)^2 &\Z/4\Z&(\Z/2\Z)^4&(\Z/3\Z)^2\\
\hline
\Omega_G&K_{12}(-2)&\Lambda_{12}(-1)&\Lambda_{14.3}(-1)&\Lambda_{15}(-1)&K_{16.3}(-1)
\end{array}
$$
%$$
%\begin{array}{l|lllll}
%G&\Z/3\Z&\Z/4\Z&(\Z/3\Z)^2&(\Z/2\Z)^2&(\Z/2\Z)^4\\
%\hline
%\Omega_G&K_{12}(-2)&\Lambda_{14.3}(-1)&K_{16.3}(-1)&\Lambda_{12}(-1)&\Lambda_{15}(-1)
%\end{array}
%$$
where $K_{12}(-2)$ is the Coxeter-Todd lattice, the lattices
$\Lambda_n$ are laminated lattices and the lattice $K_{16.3}$ is a special
 sublattice of the Leech-lattice (cf. \cite{Conway Sloane} and
\cite{Plesken Pohst} for a description). All these give lattice packings
which are very dense.\\
In each case, one can compute the invariant lattice and its
orthogonal complement by using an elliptic K3 surface. In the case
of $G=\Z/2\Z\times \Z/4\Z,(\Z/2\Z)^2,\Z/4\Z$ we use always the
same elliptic fibration, with six fibers of type $I_4$ and with
symplectic automorphisms group isomorphic to $\Z/4\Z\times \Z/4\Z$
generated by sections, then we consider its subgroups. The K3
surface admitting this elliptic fibration is the Kummer surface
${\rm Km}(E_{\sqrt{-1}}\times E_{\sqrt{-1}})$ and the group of its
automorphisms is described in \cite{Keum Kondo: automorphisms of
product elliptic curves}. This fibration admits also symplectic
automorphisms group isomorphic to $(\Z/2\Z)^i$, $i=3,4$, in this
case some of the automorphisms come from automorphisms of the base
of the fibration, i.e. of $\PP^1$. Hence also for these two groups
we are able to
compute the invariant sublattice and its orthogonal complement.\\
Once we know the lattices $\Omega_G$ we can describe all the
possible N\'eron--Severi groups of algebraic K3 surfaces with
minimal Picard number and group of symplectic automorphisms $G$.
Moreover we can prove that if there exists a K3 surface with one
of those lattices as N\'eron--Severi group, then it admits a
certain finite abelian group $G$ as group of symplectic
automorphisms. These two facts are important for the
classification of K3 surfaces with symplectic automorphisms group
(cf. \cite{bert Nikulin involutions}, \cite{symplectic prime}), in
fact they give information on the coarse moduli space of algebraic
K3 surfaces with symplectic automorphisms. For example an
algebraic K3 surface $X$ has $G$ as a symplectic automorphism
group if and only if
$\Omega_G\subset NS(X)$ (cf. \cite[Theorem 4.15]{Nikulin symplectic}) and moreover $\rho(X)\geq 1+\rank(\Omega_G)$.\\
The paper is organized as follows: in the Sections \ref{basic} and
\ref{ellipticfi} we recall some basic results on K3 surfaces and
elliptic fibrations, in Section \ref{specialauto} we show how to
find elliptic K3 surfaces with symplectic automorphisms groups
$(\Z/2\Z)^i$, $i=3,4$. Section \ref{fibrkondo} recalls some
facts on the elliptic fibration described by Keum and Kondo in
\cite{Keum Kondo: automorphisms of product elliptic curves} and
contains a description of the lattices $H^2(X,\Z)^G$ and
$(H^2(X,\Z)^G)^{\perp}$ in the cases
$G=\Z/4\Z\times\Z/4\Z,\Z/2\Z\times \Z/4\Z,\Z/4\Z,$$(\Z/2\Z)^i$,
$i=2,3,4$. In Section \ref{othercases} we give the equations
of the elliptic fibrations for the remaining $G$ and compute the
invariant lattice and its orthogonal complement. Finally
Section \ref{families} describes the N\'eron--Severi group of K3
surfaces with finite abelian symplectic automorphism group and
deals with the moduli spaces. In the Appendix we describe briefly
the elliptic fibrations which can be used to compute the lattices
$\Omega_G$ for the group $G$ which are not analyzed in Section
\ref{fibrkondo} and we give a basis for these lattices. The proof
of these results is essentially the same as the proof of
Proposition \ref{NS and T case 4 4} in Section \ref{fibrkondo}.\\

\section{Basic results}\label{basic}

Let $X$ be a K3 surface. The second cohomology group of $X$ with
integer coefficients, $H^2(X,\Z)$, with the pairing induced by the
cup product is a lattice isometric to
$\Lambda_{K3}:=E_8(-1)^2\oplus U^3$ (the K3 lattice), where $U$ is
the hyperbolic rank 2 lattice with pairing
$\left[\begin{array}{rr}0&1\\1&0\end{array}\right]$ and $E_8(-1)$
is the rank 8 negative definite lattice associated to the Dynkin diagram $E_8$ (cf.
\cite{bpv}). Let $g$ be an automorphism of $X$. It induces an
action $g^*$, on $H^2(X,\Z)$, which is an isometry. This induces
an isometry on $H^2(X,\C)$ which preserves the Hodge
decomposition. In particular $g^*(H^{2,0}(X))=H^{2,0}(X)$.
\begin{defi}
An automorphism $g\in {\Aut}(X)$ is \textbf{symplectic} if and
only if $g^*_{|H^{2,0}(X)}=Id_{|H^{2,0}(X)}$ (i.e.
$g^*(\omega_X)=\omega_X$ with $H^{2,0}(X)=\mathbb{C}\omega_X$). We
will say that a group of automorphisms \textbf{acts
symplectically} on $X$ if all the elements of the group are
symplectic automorphisms.
\end{defi}
In the following we will denote by $NS(X)$ the N\'eron--Severi
group of $X$, which is $H^{1,1}(X)\cap H^2(X,\Z)$, and by $T_X$
the transcendental lattice, which is the orthogonal complement of
$NS(X)$ in $H^2(X,\Z)$.
\begin{rem}\label{rem: Tx fixed by isometries}{\rm \cite[Theorem 3.1 b)]{Nikulin symplectic} An automorphism $g$ of $X$ is symplectic if and only if
$g^*_{|T_X}=Id_{T_X}$.\erem}\end{rem} In \cite{Nikulin symplectic}
the finite abelian groups acting symplectically on a K3 surface
are listed and many properties of this action are given. Here we
recall the most important. Let $G$ be a finite abelian group
acting symplectically on a K3 surface $X$, then
\begin{itemize}\item $G$ is one of the following fourteen groups:
\begin{eqnarray}\label{formula: symplectic group}\begin{array}{l}\Z/n\Z,\ 2\leq n\leq 8,\
(\Z/m\Z)^2,\ m=2,3,4,\\ \Z/2\Z\times \Z/4\Z,\ \Z/2\Z\times
\Z/6\Z,\ (\Z/2\Z)^i,\ i=3,4;\end{array}\end{eqnarray}
\item the
desingularization of the quotient $X/G$ is a K3 surface;
\item the
action induced by the automorphisms of $G$ on $H^2(X,\Z)$ is
unique up to isometries. In particular the lattices $H^2(X,\Z)^G$
and $(H^2(X,\Z)^G)^{\perp}$ depend on $G$ but they do not depend on the surface
X (up to isometry).
\end{itemize}
The last property implies, that we can consider a particular K3
surface $X$ admitting a finite abelian group $G$ of symplectic
automorphisms to analyze the isometries induced on
$H^2(X,\Z)\simeq \Lambda_{K3}$ and to describe the lattices
$H^2(X,\Z)^G$ and $\Omega_G=(H^2(X,\Z)^G)^{\perp}$. In particular
since $T_X$ is invariant under the action of $G$ (by Remark
\ref{rem: Tx fixed by isometries}),
$$H^2(X,\Z)^G\hookleftarrow NS(X)^G\oplus T_X$$ and the inclusion
has finite index, so considering the orthogonal lattices, we
obtain
$$(H^2(X,\Z)^G)^{\perp}=(NS(X)^G)^{\perp}.$$For a systematical approach to the problem,
how to construct K3 surfaces admitting certain symplectic
automorphisms, we will consider K3 surfaces admitting an
elliptic fibration. We recall here some basic facts.\\

Let $X$ be a K3 surface admitting an elliptic fibration with a
section, i.e. there exists a morphism $X\ra\mathbb{P}^1$ such that
the generic fiber is a non singular genus one curve and such that
there exists a section $s:\mathbb{P}^1\ra X$, which we call the
zero section. There are finitely many singular fibers, which can
also be reducible. In the following we will consider only singular
fibers of type $I_n$, $n\geq 0$, $n\in\mathbb{N}$. The fibers of
type $I_1$ are curves with a node, the fibers of type $I_2$ are
reducible fibers made up of two rational curves meeting in two
distinct points, the fibers of type $I_n$, $n>2$ are made up of
$n$ rational curves meeting as a polygon with $n$ edges. We will
call $C_0$ the irreducible component of a reducible fiber which
meets the zero section. The irreducible components of a fiber of
type $I_n$ are called $C_i$, where $C_i\cdot C_{i+1}=1$ and
$i\in\Z/n\Z$. Under the assumption $C_0\cdot s=1$, these
conditions identify the components completely once the component
$C_1$ is chosen, so these conditions identify the components up to
the transformation $C_i\leftrightarrow C_{-i}$ for each
$i\in\Z/n\Z$. All the components of a reducible fiber of type
$I_n$ have multiplicity one, so a section can intersect a fiber of
type $I_n$ in any component.\\

The set of the sections of an elliptic fibration form a group (the
Mordell--Weil group), with
the group law which is induced by the one on the fibers.\\
Let $Red$ be the set $Red=\{v\in \mathbb{P}^1\,|\, \mbox{ the
fiber }F_v\mbox{ is reducible}\}$. Let $r$ be the rank of the
Mordell--Weil group (recall that if there are no sections of
infinite order then $r=0$) and let $\rho=\rho(X)$ denote the
Picard number of the surface $X$. Then
$$\rho(X)=\rank NS(X)=r+2+\sum_{v\in Red}(m_v-1)$$ (cfr. \cite[Section
7]{shioda on mordell-weil lattice}) where $m_v$ is the number of
irreducible components of the fiber $F_v$.

%The N\'eron--Severi lattice of $X$ is generated over $\Q$ by $F$,
%by the classes of irreducible components of the reducible fiber,
%which do not intersect the zero section and by the sections.
\begin{defi} The {\bf trivial lattice} $\Tr_X$ (or ${\Tr}$) of
an elliptic fibration on a surface is the lattice generated by the
class of the fiber, the class of the zero section and the classes
of the irreducible components of the reducible fibers, which do not
intersect the zero section.\end{defi} The lattice ${\Tr}$
admits $U$ as sublattice and its rank is $\rank ({\rm
Tr})=2+\sum_{v\in Red}(m_v-1)$. Recall that $NS(X)\otimes \Q$ is
generated by ${\Tr}$ and the sections of infinite order.

\begin{theorem}{\rm \cite[Theorem 1.3]{shioda on mordell-weil
lattice}} The Mordell--Weil group of the elliptic fibration on the
surface $X$ is isomorphic to the quotient $NS(X)/{\Tr}:=E(K)$.
\end{theorem}

In Section 8 of \cite{shioda on mordell-weil lattice} a pairing on
$E(K)$ is defined. The value of this pairing on a section $P$
depends only on the intersection between the section $P$ and the
reducible fibers and between $P$ and the zero section. Now we
recall the definition and the properties
of this pairing.\\
Let $E(K)_{\tors}$ be the set of the torsion elements in the group
$E(K)$.
\begin{lemma}{\rm \cite[Lemma 8.1, Lemma 8.2]{shioda on mordell-weil lattice}} For
any $P\in E(K)$ there exists a unique element $\phi(P)$ in
$NS(X)\otimes
\Q$ such that: \\
i) $\phi(P)\equiv (P) \mod ({\Tr}\otimes \Q)$ (where $(P)$ is the class of $P \mod ({\Tr}\otimes \Q)$)\\
ii) $\phi(P)\perp {\Tr}$.\\
The map $\phi:E(K)\ra NS(X)\otimes \Q$ defined above is a group
homomorphism such that $Ker(\phi)=E(K)_{tor}.$
\end{lemma}
\begin{theorem}\label{theorem: height formula}{\rm \cite[Theorem 8.4]{shioda on mordell-weil
lattice}} For any $P,Q\in E(K)$ let $\langle
P,Q\rangle=-\phi(P)\cdot\phi(Q)$ (where $\cdot$ is induced on
$NS(X)\otimes\Q$ by the cup product). This defines a symmetric
bilinear pairing on $E(K)$, which induces the structure of a
positive definite lattice on
$E(K)/E(K)_{tor}$.\\
In particular if $P\in E(K)$, then $P$ is a torsion section if and
only if $\langle P,
P\rangle=0$.\\
For any $P,Q\in E(K)$ the pairing $\langle-,-\rangle$ is
$$
\begin{array}{lll}
\langle P, Q \rangle&=&\chi+P\cdot s+Q\cdot s-P\cdot Q-\sum_{v\in
Red}{\rm contr}_v(P,Q)\\
\langle P, P \rangle&=&2\chi+2(P\cdot s)-\sum_{v\in Red}{\rm
contr}_v(P)
\end{array}$$
 where $\chi$ is the holomorphic Euler characteristic of the surface and the rational numbers ${\rm contr}_v(P)$ and ${\rm contr}_v(P,Q)$ are given in the table below
\begin{eqnarray}
\begin{array}{llllll}
&I_2&I_n&I_n^*&IV^*&III^*\\
{\rm contr}_v(P)&2/3&i(n-i)/n &\left\{\begin{array}{l} 1\mbox{ \rm if
}
i=1\\1+n/4\mbox{ \rm if } i=n-1\mbox{ \rm or }i=n\end{array}\right.&4/3&3/2\\
{\rm contr}_v(P,Q)&1/3&i(n-j)/n &\left\{\begin{array}{l} 1/2\mbox{ \rm
if } i=1\\2+n/4\mbox{ \rm if } i=n-1\mbox{\rm  or
}i=n\end{array}\right.&2/3&-
\end{array}
\end{eqnarray}
where the numbering of the fibers is the one described before, $P$
and $Q$ meet the fiber in the component $C_i$ and $C_j$ and $i\leq
j$.\end{theorem} The pairing defined in the theorem is called the
\textbf{height pairing}. This pairing will be used to determine
the intersection of the torsion sections of the elliptic
fibrations with the irreducible components of the reducible
fibers.

\section{Elliptic fibrations and symplectic automorphisms}\label{ellipticfi}

Using elliptic fibrations one can describe the action of a
symplectic automorphism induced by a torsion section on the
N\'eron--Severi group. Since a symplectic automorphism acts as the
identity on the transcendental lattice, we can describe the action
of the symplectic automorphism on $NS(X)\oplus T_X$, which is a
sublattice of finite index of $H^2(X,\Z)$. Moreover, knowing the
discriminant form of $NS(X)$ and of $T_X$, one can explicitly find
a basis for the lattice $H^2(X,\Z)$, and so one can describe the
action of the symplectic automorphism on the
lattice $H^2(X,\Z)$.\\
Let $X$ be a K3 surface admitting an elliptic fibration with a
section, then the N\'eron--Severi group of $X$ contains the
classes $F$ and $s$, which are respectively the class of the fiber
and of the section. Let us suppose that $t$ is a section of an
elliptic K3 surface and let $\sigma_t$ be the automorphism,
induced by $t$, which fixes the base of the fibration (so it fixes
each fiber) and acts on each fiber as translation by $t$. This
automorphism is symplectic and if the section $t$ is an
$n$-torsion section (with respect to the group law of the elliptic
fibration) the automorphism $\sigma_t$ has order $n$. More in
general let us consider a K3 surface $X$ admitting an elliptic
fibration $\mathcal{E}_X$ with torsion part of the Mordell--Weil
group equal to a certain abelian group $G$. Then the elements of
 $\tors(\MW(\mathcal{E}_X))$ induce symplectic automorphisms which
commute. So we obtain the following:

\begin{lemma}\label{rem: if G=torsMW the G is symplectic on X}If $X$ is a K3 surface with an elliptic fibration
$\mathcal{E}_X$ and $\tors(\MW(\mathcal{E}_X))=G$, then $X$ admits
$G$ as abelian group of symplectic automorphisms and these
automorphisms are induced by the torsion sections of
$\mathcal{E}_X$.\end{lemma} Now we want to analyze the action of
the symplectic automorphisms induced by torsion sections on the
classes generating the N\'eron--Severi group of the elliptic
fibration. By definition, $\sigma_t$ fixes the class $F$. Since it
acts as a translation on each fiber, it sends on each fiber the
intersection of the fiber with the zero section, to the
intersection of the fiber with the section $t$. Globally
$\sigma_t$ sends the section $s$ to the section $t$. More in
general the automorphism $\sigma_t$ sends sections to sections and
the section $r$ will be send to $r+t-s$, where $+$ and $-$ are
the operations in the Mordell--Weil group. To complete the
description of the action of $\sigma_t$ on the N\'eron--Severi
group we have to describe its action on the reducible fibers. The
automorphism $\sigma_t$ restricted to a reducible fiber has to be
an automorphism of it. This fact imposes restrictions on the
existence of a torsion section and of certain reducible fibers in
the same elliptic fibration.
\begin{lemma}\label{rem: action of sigmat non-trivial other
C0} Let $t$ be an $n$-torsion section and $F_v$ be a reducible fiber
of type $I_d$ of the fibration. Let us suppose that the section
$t$ meets the fiber $F_v$ in the component $C_i$.\\ Let $r$ be the
minimal number such that $(\sigma_t)_{|F_v}^r(C_j)=C_j$, $\forall
j\in\Z/d\Z$.\\
If $\gcd(d,n)=1$ then $r=1$, $i=0$ and $\sigma_t(C_j)=C_j$, $\forall j\in\Z/d\Z$.\\
If $\gcd(d,n)\neq 1$ and $i=0$, then $r=1$ and $\sigma_t(C_j)=C_j$, $\forall j\in\Z/d\Z$.\\
If $\gcd(d,n)\neq 1$ and $i\neq 0$, then $r=d/\gcd(d,i)$,
$r|\gcd(d,n)$ and $\sigma_t(C_j)=C_{j+i}$.\end{lemma} \bprf We
prove first that if $t\cdot C_i=1$, then $\sigma_t(C_j)=C_{i+j}$.
The section $s$ meets the reducible fiber in the component $C_0$,
so $1=s\cdot C_0=\sigma_t(s)\cdot \sigma_t(C_0)=t\cdot
\sigma_t(C_0)$. But $\sigma_t(C_0)$ has to be a component of the
same fiber (because $\sigma_t$ fixes the fibers) and has to be the
component with a non-trivial intersection with $t$, so
$\sigma_t(C_0)=C_i$. Moreover $1=C_0\cdot C_{1}=\sigma_t(C_0)\cdot
\sigma_t(C_1)=C_i\cdot \sigma_t(C_{1})$ and fixing an orientation
on the fiber we obtain $\sigma_t(C_1)=C_{i+1}$. Hence, more in
general, one has $\sigma_t(C_j)=C_{j+i}$. Moreover if
$(\sigma_t)_{|F_v}(C_j)=C_j$ for each $j\in\Z/d\Z$, then $t\cdot C_0=1$.\\
Let $r$ be the minimal number such that
$(\sigma_t)^r_{|F_v}(C_j)=C_j$ for each $j\in\Z/d\Z$ . Since the
index of the components of a reducible fibers of type $I_n$ are
defined modulo $d$, it is clear that $\sigma_t^d(C_j)=C_j$ for
each $C_j$ component of $F_v$. Observe that $r\leq d$ and in
particular $r|d$. Now since $\sigma_t$ has order $n$
 and $\sigma_t^r(C_j)=C_j$ for each $ j$, we obtain $r|n$.
Hence $r|\gcd(d,n)$. If $\gcd(d,n)=1$, then $r=1$,
$\sigma_t(C_j)=C_j$ for each $j$ and so $i=0$. If $\gcd(d,n)\neq
1$, then $r$ has to be the smallest positive number such that
$ri\equiv 0\mod d$, this implies $r=d/\gcd(d,i)$.\eprf

In the following we always use the notation: if $t_1$ is an
$n$-torsion section, then $t_h$ corresponds to the sum, in the
Mordell--Weil group law, of $h$ times $t_1$. Moreover $C_i^{(j)}$
is the $i$-th component of the $j$-th reducible fiber.

\section{The cases $G=(\Z/2\Z)^3$ and $G=(\Z/2\Z)^4$}\label{specialauto}
We have seen (cf. Lemma \ref{rem: if G=torsMW the G is symplectic
on X}) that an example of a K3 surface with $G$ as group of
symplectic automorphisms is given by an elliptic fibration with
$G$ as torsion part of the Mordell--Weil group. There are twelve
groups which appear as torsion part of the Mordell--Weil group of
an elliptic fibration on a K3 surface. In particular the groups
$G=(\Z/2\Z)^3$ and $G=(\Z/2\Z)^4$ are groups acting symplectically
on a K3 surface, but they can not be realized as the torsion part
of the Mordell--Weil group of an elliptic fibration on a K3
surface (cf. \cite{shimada}). Hence to find examples of K3
surfaces with one of these groups as group of symplectic
automorphisms, we have to use a different construction. One
possibility is to consider a K3 surface with an elliptic fibration
$\mathcal{E}_X$ with $\MW(\mathcal{E}_X)=(\Z/2\Z)^2$ and to find
one (resp. two) other symplectic involutions which commute with
the ones induced by torsion sections.

\subsection{The group $G=(\Z/2\Z)^{3}$ acting symplectically on
an elliptic fibration}

We start considering an elliptic K3 surface with two 2-torsion
sections. An equation of such an elliptic fibration is given by

\vspace{-0.5cm}

\begin{equation}\label{equation 2 torsion} y^2=x(x-p(\tau))(x-q(\tau)),\ \ \deg(p(\tau))=\deg(q(\tau))=4,\
\tau\in\C.\end{equation}

\vspace{-0.1cm}

Then we consider an involution on the base of the fibration
$\mathbb{P}^1$, which preserves the fibration. This involution
fixes two points of the basis. Up to the choice of the coordinates
on $\mathbb{P}^1$, we can suppose that the involution on the basis
is $\sigma_{\mathbb{P}^1,a}:\tau\mapsto -\tau$. So we consider on
the K3 surface the involution $(\tau,x,y)\mapsto (-\tau,x,-y)$.
Since this map has to be an involution of the surface, it has to
fix the equation of the elliptic fibration. Moreover the
involution $\sigma_{\mathbb{P}^1,a}$ has to commute with the
involutions induced by the torsion sections. This implies that it
has to fix the curves corresponding to the torsion sections
$t:\tau\mapsto(p(\tau),0)$ and $u:\tau\mapsto (q(\tau),0)$. The
equation of such an elliptic fibration  is

\vspace{-0.5cm}

\begin{eqnarray}\label{equation 2 torsion 2 torsion+involution}\begin{array}{c}y^2=x(x-p(\tau))(x-q(\tau))\ \ \mbox{ with}\\\deg(p(\tau))=\deg(q(\tau))=4\mbox{ and } p(\tau)=p(-\tau),\ \ q(\tau)=q(-\tau),\end{array}\end{eqnarray}

%\vspace{-0.2cm}

i.e.
 $y^2=x(x-(p_4\tau^4+p_2\tau^2+p_0))(x-(q_4\tau^4+q_2\tau^2+q_0)),\
\ \ p_4,p_2,p_0,q_4,q_2,q_0\in\C.$\\ The involution
$\sigma_{\mathbb{P}^1,a}$ fixes the four curves corresponding to
the sections in the torsion part of the Mordell--Weil group and
the two fibers over the points $0$ and $\infty$ of $\mathbb{P}^1$.
On these two fibers the automorphism is not the identity (because
it sends $y$ to $-y$). So it fixes only the eight intersection
points between these four sections and the two fibers. It fixes
eight isolated
points, so it is a symplectic involution.\\
To compute the moduli of this family of surfaces we have to
consider that the choice of the automorphism on $\mathbb{P}^1$
corresponds to a particular choice of the coordinates, so we can not
act on the equation with all the automorphisms of $\mathbb{P}^1$.
We choose coordinates, such that the involution $\sigma_{\mathbb{P}^1,a}$
fixes the points $(1:0)$ and
$(0:1)$ on $\mathbb{P}^1$.

The space of the automorphisms of $\mathbb{P}^1$ commuting with
$\sigma_{\mathbb{P}^1,a}$ has dimension one. Moreover we can act
on the equation \eqref{equation 2 torsion 2 torsion+involution}
with the transformation $(x,y)\mapsto(\lambda^2x,\lambda^3y)$ and
divide by $\lambda^6$. Since there are six parameters in the
equation \eqref{equation 2 torsion 2 torsion+involution}, we find
that the number of moduli of this family is $6-2=4$. Collecting
these results, the properties of the family satisfying the
equation \eqref{equation 2 torsion 2 torsion+involution} are the
following:
$$\begin{array}{c|c|c}
\mbox{discriminant}& \mbox{singular fibers}& \mbox{moduli}\\
\hline
p(\tau)^2q(\tau)^2(p(\tau)-q(\tau))^2&12\,I_2&4
\end{array}
$$
Since the number of moduli of this family is four, the Picard
number of the generic surface is $\rho\leq 16$.\\
The automorphisms induced by the translation by the two 2-torsion
sections and by $\sigma_{\mathbb{P}^1,a}$ fix the fiber and the
sum of the 2-torsion sections with the zero section, so the
invariant sublattice of the N\'eron--Severi group contains at
least two classes. Its orthogonal complement in the
N\'eron--Severi group is a lattice which does not depend on the
surface (as explained before) and which has rank 14 \cite{Nikulin
symplectic}. Hence the Picard number of the surfaces in this
family is at least 2+14=16.\\ This implies that the generic
member of the family has Picard number 16. The trivial lattice of
this fibration has rank 14, so there are two linearly independent
sections of infinite order which generate the free part of the
Mordell--Weil group.

\subsection{The group $G=(\Z/2\Z)^{4}$ acting symplectically on
an elliptic fibration}

As in the previous section we construct an involution which
commutes with the three symplectic involutions $\sigma_t$,
$\sigma_u$ (the involutions induced by the torsion sections $t$ and
$u$) and $\sigma_{\mathbb{P}^1,a}$ of the surfaces described by
the equation \eqref{equation 2 torsion 2 torsion+involution}. So
we consider two commuting involutions on $\mathbb{P}^1$ which
commute with the involutions induced by the 2-torsion sections. Up
to the choice of the coordinates of $\mathbb{P}^1$ we can suppose
that the involutions on $\mathbb{P}^1$ are $\tau\mapsto -\tau$ and
$\tau\mapsto 1/\tau$. We call the corresponding involutions on the
surface $\sigma_{\mathbb{P}^1,a}:(\tau,x,y)\mapsto (-\tau,x,-y)$
and $\sigma_{\mathbb{P}^1,b}:(\tau,x,y)\mapsto (1/\tau,x,-y)$. As
before requiring that these involutions commute with the
involutions induced by the torsion sections means that each of
these involutions fixes the torsion sections. Elliptic K3 surfaces
with the properties described have the following
equation:

\vspace{-0.7cm}

\begin{eqnarray}\label{equation 2 torsion 2 torsion+2
involutions}\begin{array}{c}y^2=x(x-p(\tau))(x-q(\tau))\ \ \mbox{
with}\\ \deg(p(\tau))=\deg(q(\tau))=4,\ \
p(\tau)=p(-\tau)=p(\frac{1}{\tau}),\ \
q(\tau)=q(-\tau)=q(\frac{1}{\tau}).\end{array}\end{eqnarray} i.e.
$y^2=x(x-(p_4\tau^4+p_2\tau^2+p_4))(x-(q_4\tau^4+q_2\tau^2+q_4)),\
\ p_4,p_2,q_4,q_2\in\C.$\\ As before the choice of the involutions
of $\mathbb{P}^1$ implies that the admissible transformations on
the equation have to commute with $\sigma_{\mathbb{P}^1,a}$ and
$\sigma_{\mathbb{P}^1,b}$ . The only possible transformation is
the identity so
%preserve the points $(1:0)$, $(0:1)$,
%$(1:1)$, $(1:-1)\in\mathbb{P}^1$. So there are
no transformations
of $\mathbb{P}^1$ can be applied to the equation
\eqref{equation 2 torsion 2 torsion+2 involutions}. The only
admissible transformation on that equation is
$(x,y)\mapsto(\lambda^2x,\lambda^3y)$.\\
Collecting these results, the properties of the family satisfying
the equation \eqref{equation 2 torsion 2 torsion+2 involutions}
are the following:

%\vspace{-0.5cm}

$$\begin{array}{c|c|c}
\mbox{discriminant}& \mbox{singular fibers}& \mbox{moduli}\\
\hline
p(\tau)^2q(\tau)^2(p(\tau)-q(\tau))^2&12\,I_2&3
\end{array}
$$
As in the previous case the comparison between the number of
moduli of the family and the rank of the trivial lattice implies
that the free part of the Mordell--Weil group is $\Z^3$.

\section{An elliptic fibration with six fibers of type $I_4$}\label{fibrkondo}

An equation of an elliptic K3 surface with six fibers of type
$I_4$ is the following
\begin{eqnarray}\label{formula: Weierstrass equation 4 4}
y^2=x(x-\tau^2\sigma^2)\left(x-\frac{(\tau^2+\sigma^2)^2}{4}\right).\end{eqnarray}
This equation is well known, for example it is described in
\cite[Section 2.3.1]{Top Yui}, where it is constructed considering
the K3 surface as double cover of a rational elliptic surface with
two fibers of type $I_2$ and two fibers of type $I_4$. One can
find this equation also considering the equation of an elliptic K3
surface with two 2-torsion sections (i.e the equation
\eqref{equation 2 torsion}) and requiring that the tangent lines
to the elliptic curve defined by the equation \eqref{equation 2
torsion} in the two rational points of the elliptic surface pass
respectively through the two points of order two $(p(\tau),0)$ and
$(q(\tau),0)$.

\begin{figure}[t]
\begin{center}
\includegraphics[scale=1.25]{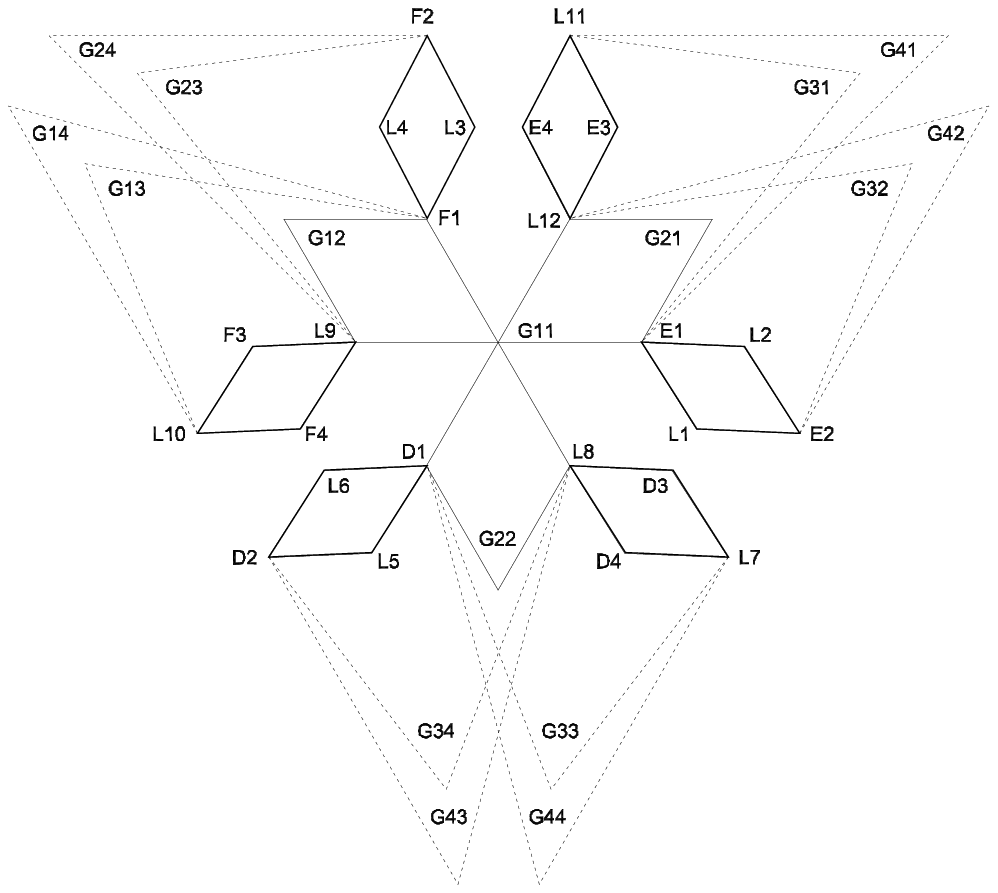}
\caption{}
\label{ali3}
\end{center}
\end{figure}

We will call $X_{(\Z/4\Z)^2}$ the elliptic K3 surface described by
the equation \eqref{formula: Weierstrass equation 4 4}. It has two
4-torsion sections $t_1$ and $u_1$, which induce two commuting
symplectic automorphisms $\sigma_{t_1}$ and $\sigma_{u_1}$. We
will analyze the properties of this surface to study the group of
symplectic automorphisms $G=(\Z/4\Z)^2$ and also its subgroups
$\Z/2\Z\times \Z/4\Z=\langle\sigma^2_{t_1}, \sigma_{u_1}\rangle$,
$(\Z/2\Z)^2=\langle
\sigma_{t_1}^2,\sigma_{s_1}^2\rangle$, $\Z/4\Z=\langle \sigma_{t_1}\rangle$, $\Z/2\Z=\langle\sigma_{t_1}^2\rangle$ which act symplectically too.\\
It is more surprising that the equation \eqref{formula:
Weierstrass equation 4 4} appears as a specialization also of the
surfaces described in \eqref{equation 2 torsion 2
torsion+involution} and \eqref{equation 2 torsion 2 torsion+2
involutions}. So it admits also the group $G=(\Z/2\Z)^4$ as group
of symplectic automorphisms. The automorphisms
$\sigma_{\mathbb{P}^1,a}:(\tau,x,y)\mapsto (-\tau,x,-y)$ and
$\sigma_{\mathbb{P}^1,b}:(\tau,x,y)\mapsto (1/\tau,x,-y)$ commute
with the automorphism induced by the two 2-torsion sections $t_2$
and $u_2$.\\ The automorphisms group ${\Aut}(X_{(\Z/4\Z)^2})$
of the surface $X_{(\Z/4\Z)^2}$ is described in \cite{Keum Kondo:
automorphisms of product elliptic curves}. Keum and Kondo prove
that $X_{(\Z/4\Z)^2}$ is the Kummer surface ${\Km}(E_{\sqrt{-1}}\times E_{\sqrt{-1}})$. They consider forty
rational curves on the surface: sixteen of them, (the ones called
$G_{i,j}$), form a Kummer lattice, and the twentyfour curves
$G_{i,j}$, $E_i$, $i=1,2,3,4$, $F_j$, $j=1,2,3,4$ form the so
called {\it double Kummer}, a lattice which is a sublattice of the
N\'eron--Severi group of ${\Km}(E_1\times E_2)$ for each couple
of elliptic curves $E_1$, $E_2$. These forty curves generate
$NS({\Km}(E_{\sqrt{-1}}\times E_{\sqrt{-1}}))$ and describe
five different elliptic fibrations on $X_{(\Z/4\Z)^2}$ which have
six fibers of type $I_4$ each. A complete list of the 63 elliptic
fibrations on the surface
$X_{(\Z/4\Z)^2}$ can be found in \cite{Nishiama}.\\
In the Figure \ref{ali3} some of the intersections between the
forty curves introduced in \cite{Keum Kondo: automorphisms of
product elliptic curves} are shown.
The five elliptic fibrations are associated to the classes
(\cite[Proof of Lemma 3.5]{Keum Kondo: automorphisms of product
elliptic curves})
$$E_1+E_2+L_1+L_2,\ \ G_{11}+G_{12}+F_1+L_9,\ \ G_{11}+G_{21}+E_1+L_{12},$$
$$ G_{11}+G_{22}+D_1+L_8,\ \ G_{33}+G_{44}+D_1+L_7.$$
We identify the fibration that we will consider with one of these,
say the one with fiber $E_1+E_2+L_1+L_2$. Then we can define
$s=G_{11}$, $t=G_{13}$, $u=G_{32}$, and so $C_{0}^{(1)}=E_1$,
$C_{0}^{(2)}=L_{12}$, $C_{0}^{(3)}=L_8$, $C_{0}^{(4)}=D_1$,
$C_{0}^{(5)}=E_1$, $C_{0}^{(6)}=L_9$.\\
We will consider the group $G=(\Z/2\Z)^4$ generated by two
involutions induced by torsion sections and two involutions
induced by involutions on $\mathbb{P}^1$ (with respect to one of
these elliptic fibrations). However one can choose the last two
involutions as induced by the two torsion sections of different
elliptic fibrations (cf. Remark \ref{rem: automorphisms P1 are
induced by torsion sections}).
%However it can be
%considered also generated by the two torsion sections of
%different elliptic fibrations.

Up to now for simplicity we put $X_{(\Z/4\Z)^2}=X$.

\begin{prop}\label{NS and T case 4 4} Let $X$ be the
elliptic K3 surface with equation \eqref{formula: Weierstrass
equation 4 4}. Let $t_1$ and $u_1$ be two 4-torsion sections of
the fibration. Then $t_1\cdot C_1^{(j)}=1$ if $j=1,2,3,4$,
$t_1\cdot C_2^{(5)}=t_1\cdot C_0^{(6)}=1$, $u_1\cdot C_1^{(h)}=1$
if $h=4,5,6$, $u_1\cdot C_3^{(3)}=1$, $u_1\cdot
C_2^{(1)}=u_1\cdot C_0^{(2)}=1$.\\
A $\mathbb{Z}$-basis for the lattice $NS(X)$ is given by $F,\ s,\
t_1,\ u_1,\ C_i^{(j)},\ j=1,\ldots,6$, $i=1,2$ and $C_3^{(j)}$, $j=2,\ldots,5$.\\
The trivial lattice of the fibration is $U\oplus A_3^{\oplus 6}$.
It has index $4^2$ in the N\'eron--Severi group of $X$. The
lattice $NS(X)$ has discriminant $-4^2$ and its discriminant form
is $\Z/4\Z(-\frac{1}{4})\oplus \Z/4\Z(-\frac{1}{4})$.\\ The
transcendental lattice is
$T_{X}=\left[\begin{array}{ll}4&0\\0&4\end{array}\right]$, and has
a unique primitive embedding in the lattice
$\Lambda_{K3}.$\end{prop} \bprf The singular fibers of this
fibration are six fibers of type $I_4$ (the classification of the
type of the singular fibers is determined by the zero-locus of the
discriminant of the equation of the surface and can be found in
\cite[IV.3.1]{miranda elliptic pisa}). Hence the trivial lattice
of this elliptic fibration is $U\oplus A_3^{\oplus 6}$. Since it
has rank 20, which is the maximal Picard number of a K3 surface,
there are no sections of infinite order on this elliptic
fibration. The torsion part of the Mordell--Weil group is
$(\Z/4\Z)^2$, generated by two 4-torsion sections. We will call
these sections $t_1$ and $u_1$. By the height formula (cf. Theorem
\ref{theorem: height formula}) the intersection between a four
torsion section and the six fibers of type $I_4$ has to be of the
following type: the section has to meet four fibers in the
component $C_{i}$, with $i$ odd, a fiber in the component $C_{2}$
and a fiber in the component $C_0$. After a suitable numbering of
the fibers we can suppose that the section $t_1$ has the
intersection described in the statement. The section $u_1$ and the
section $v_1$, which corresponds to $t_1+u_1$ in the Mordell--Weil
group, must intersect the fibers in a similar way (four in the
component $C_i$ with an odd $i$, one in $C_2$ and one in $C_0$).
If $t_1\cdot C_i^{(j)}=1$ and $u_1\cdot C_h^{(j)}=1$, then $v_1$
meets the fiber in the component $C_{i+h}^{(j)}$ (it is a
consequence of the group law on the fibers of type $I_n$). The
conditions on the intersection properties of $u_1$ and $v_1$ imply
that $u_1\cdot C_1^{(h)}=1$ if $h=3,4,5,6$, $u_1\cdot
C_2^{(1)}=u_1\cdot C_0^{(2)}=0$, and hence $v_1\cdot C_1^{(j)}=1$,
$j=2,6$, $v_1\cdot C_2^{(4)}=1$, $v_1\cdot
C_3^{(h)}=1$, $h=1,5$, $v_1\cdot C_0^{(3)}=1$.\\
The torsion sections $t_1$ and $u_1$ can be written as a linear
combination of the classes in the trivial lattice with
coefficients in $\frac{1}{4}\Z$ (because they are 4-torsion
sections). Hence the trivial lattice has index $4^2$ in the
N\'eron--Severi group and so $\discr(NS(X))=4^6/4^4=4^2$. In particular
$$\begin{array}{l}t_1=2F+s-\frac{1}{4}\left(\sum_{i=1}^4(3C_1^{(i)}+2C_2^{(i)}+C_3^{(i)})+2C_1^{(5)}+4C_2^{(5)}+2C_3^{(5)}\right),\\
u_1=2F+s-\frac{1}{4}\left(\sum_{i=4}^6(3C_1^{(i)}+2C_2^{(i)}+C_3^{(i)})+C_1^{(3)}+2C_2^{(3)}+3C_3^{(3)}+2C_1^{(1)}+4C_2^{(1)}+2C_3^{(1)}\right).\end{array}
$$
It is now clear that $F,\ s,\ t_1,\ u_1,\ C_i^{(j)},\
j=1,\ldots,6$, $i=1,2$ and $C_3^{(j)}$, $j=2,\ldots,5$ is a
$\Q$-basis for the N\'eron--Severi group and since the determinant
of the intersection matrix of this basis is $4^2$, it is in fact a
$\Z$-basis. The discriminant form of the N\'eron--Severi lattice
is generated by
\begin{eqnarray}\label{formula: generators discriminant form}\begin{array}{l}d_1=\frac{1}{4}(C_1^{(1)}+2C_2^{(1)}+3C_3^{(1)}-C_1^{(3)}-2C_2^{(3)}-3C_3^{(3)}+C_1^{(6)}+2C_2^{(6)}+3C_3^{(6)}),\\
d_2=\frac{1}{4}(C_1^{(2)}+2C_2^{(2)}+3C_3^{(2)}+C_1^{(4)}+2C_2^{(4)}+3C_3^{(4)}-C_1^{(5)}-2C_2^{(5)}-3C_3^{(5)}),\end{array}\end{eqnarray}
hence the discriminant form of $NS(X)$ is
$\Z/4\Z(-\frac{1}{4})\oplus\Z/4\Z(-\frac{1}{4})$. The
transcendental lattice has to be a rank 2 positive definite
lattice with discriminant form
$\Z/4\Z(\frac{1}{4})\oplus\Z/4\Z(\frac{1}{4})$ (the opposite of
the one of the N\'eron--Severi group). This implies that
$T_{X}=\left[\begin{array}{rr}4&0\\0&4\end{array}\right].$ By
\cite[Theorem 1.14.4]{Nikulin symplectic} it admits a unique
primitive embedding in $\Lambda_{K3}$.\erem

\begin{prop}\label{omega and invariant 44} Let $G_{4,4}$ be the group
generated by $\sigma_{t_1}$ and $\sigma_{u_1}$. The invariant
sublattice of the N\'eron--Severi
group with respect to $G_{4,4}$ is isometric to $\left[\begin{array}{rr}-8&8\\8&0\end{array}\right]$.\\
Its orthogonal complement $(NS(X)^{G_{4,4}})^{\perp}$ is
$\Omega_{(\Z/4\Z)^2}:=(H^2(X,\mathbb{Z})^{G_{4,4}})^{\perp}$. It
is the negative definite lattice $\{\mathbb{Z}^{18}, Q\}$ where
$Q$ is the bilinear form obtained as the intersection form of the
classes
\begin{eqnarray*}
\begin{array}{l}
b_{1}=s-t_1,\ b_2=s-u_1,\ b_{i+2}=C_1^{(i)}-C_1^{(i+1)},\
i=1,\ldots,5,\ b_{j+7}=C_2^{(j)}-C_2^{(j+1)}, j=1,\ldots,5,\\
b_{h+11}=C_1^{(h)}-C_3^{(h)},\ h=2,\ldots,5,\
b_{17}=C_1^{(1)}-C_2^{(2)},\
b_{18}=f+s-t_1-C_1^{(1)}-C_1^{(2)}-C_2^{(2)}-C_1^{(3)}.
\end{array}
\end{eqnarray*}
The lattice $\Omega_{(\Z/4\Z)^2}$ admits a unique primitive
embedding in
the lattice $\Lambda_{K3}$.\\
The discriminant of $\Omega_{(\Z/4\Z)^2}$ is $2^{8}$ and its discriminant group is $(\Z/2\Z)^{ 2}\oplus(\Z/8\Z)^{ 2}$.\\
The group of isometries $G_{4,4}$ acts on the discriminant group
$\Omega_{(\Z/4\Z)^2}^{\vee}/\Omega_{(\Z/4\Z)^2}$ as the
identity.\\ The lattice $H^2(X,\Z)^{G_{4,4,}}$ is \vspace{-0.5cm}

$$\left[\begin{array}{rrrr}
  4& 6& 0& 0\\  6& 4& 6& 4\\  0& 6& 4& 0\\
  0& 4& 0& 0\end{array}\right]
$$
and it is an overlattice of $NS(X)^{G_{4,4}}\oplus T_{X}\simeq
\left[\begin{array}{rr}-8&8\\8&0\end{array}\right]\oplus
\left[\begin{array}{ll}4&0\\0&4\end{array}\right]$ of index
four.\end{prop} \bprf Let $v_1=t_1+u_1$, $w_1=t_2+u_1$,
$z_1=t_3+u_1$. By definition, the automorphisms induced by the
torsion sections fix the class of the fiber, so $F\in
NS(X)^{G_{4,4}}$. Moreover the group $G_{4,4}$ acts on the section
fixing the class $s+\sum_{i=1}^3
\left(t_i+u_i+v_i+w_i+z_i\right)$. The action of the group
$G_{4,4}$ is not trivial on the components of the reducible
fibers, so the lattice $\langle F, s+\sum_{i=1}^3
\left(t_i+u_i+v_i+w_i+z_i\right)\rangle$ is a sublattice of
$NS(X)^{G_{4,4}}$ of finite index. Since the orthogonal complement of a
sublattice  is always a primitive sublattice, we have $\langle F,
s+\sum_{i=1}^3
\left(t_i+u_i+v_i+w_i+z_i\right)\rangle^{\perp}=(NS(X)^{G_{4,4}})^{\perp}$.
In this way one can compute the classes generating the lattice
$(NS(X)^{G_{4,4}})^{\perp}$. A basis for this lattice is given by
the classes $b_i$, moreover the lattice $NS(X)^{G_{4,4}}$ is
isometric to $((NS(X)^{G_{4,4}})^{\perp})^{\perp}$, and a
computation shows that it is isometric to
$\left[\begin{array}{rr}-8&8\\8&0\end{array}\right]$.\\
To find the lattice $H^2(X,\Z)^{G_{4,4}}$ we consider the
orthogonal complement of
$(H^2(X),\Z)^{G_{4,4}})^{\perp}\simeq(NS(X)^{G_{4,4}})^{\perp}$
inside the lattice $H^2(X,\Z)$. Since we know the generators of
the discriminant form of $NS(X)$ and of $T_{X}$, we can construct
a basis of $H^2(X,\Z)$. Indeed let $a_1$ and $a_2$ be the
generators of $T_{X}$, then the classes $F$, $s$, $t_1$, $u_1$,
$C_i^{(j)},\ j=1,\ldots,6$, $i=1,2$, $C_3^{(j)}$, $j=2,\ldots,5$,
$a_1/4+d_1$ and $a_2/4+d_2$ ($d_i$ as in \eqref{formula:
generators discriminant form}) form a $\Z$-basis of $H^2(X,\Z)$.
The classes $b_i$, $i=1,\ldots, 18$ generate
$(H^2(X,\Z)^{G_{4,4}})^{\perp}$ and are expressed as a linear
combination of the $\Z$-basis of $H^2(X,\Z)$ described above. The
orthogonal complement of these classes in $H^2(X,\Z)$ is the
lattice $H^2(X,\Z)^{G_{4,4}}$. Since we found the lattice
$\Omega_{(\Z/4\Z)^2}$, by a computer computation it is easy  to
find its discriminant group which is $(\Z/2\Z)^{
2}\oplus(\Z/8\Z)^{ 2}$ and to study the action of $G_{4,4}$ on it,
which
 is trivial.

%A computation
%with the computer shows that the action of $G_{4,4}$ is trivial on
%the discriminant group
%$\Omega_{(\Z/4\Z)^2}^{\vee}/\Omega_{(\Z/4\Z)^2}=(\Z/2\Z)^{
%2}\oplus(\Z/8\Z)^{ 2}$.
%The discriminant
%group $\Omega_{(\Z/4\Z)^2}^{\vee}/\Omega_{(\Z/4\Z)^2}=(\Z/2\Z)^{
%2}\oplus(\Z/8\Z)^{ 2}$ is generated by ({\bf hai i generatori?})
%and a computation shows that the action of $G_{4,4}$ is trivial.
Moreover, the lattice $\Omega_{(\Z/4\Z)^2}$ satisfies the
hypothesis of \cite[Theorem 1.14.4]{Nikulin bilinear}, so  it
admits a unique primitive embedding in $\Lambda_{K3}$. \erem

\begin{prop}\label{omega and invariant 24 22 4}{\bf 1)} Let $G_{2,4}$ be the group generated by $\sigma_{t_2}$ and
$\sigma_{u_1}$. The invariant sublattice of the N\'eron--Severi
group with respect to $G_{2,4}$ is isometric to $U(4)\oplus
\left[\begin{array}{rr}-4&0\\0&-4\end{array}\right].$
%su mathematica ho la base... ma non so se serve.
Its orthogonal complement $(NS(X)^{G_{2,4}})^{\perp}$ is
$\Omega_{\Z/2\Z\times\Z/4\Z}:=(H^2(X,\mathbb{Z})^{G_{2,4}})^{\perp}$.
It is the negative definite lattice $\{\mathbb{Z}^{16}, Q\}$ where
$Q$ is the bilinear form obtained as the intersection form of the
classes
\begin{eqnarray*}
\begin{array}{ll}
b_1=s-u_1,\ \ \ \ b_{i}=C_3^{(i)}-C_1^{(i)},\!\!\!\!&i=2,\ldots,
5,\
\ \ \ \ \ \ b_{j+3}=C_1^{(j)}-C_1^{(j+1)},\ j=3,4,5,\\
b_{h+6}=C_2^{(h)}-C^{(h)}_2,\
h=3,4,5,\!\!\!\!&b_{12}=C_2^4-C_1^3,\ \
b_{13}=C_0^{(3)}-C_1^{(4)},\ \
b_{14}=C_1^{(1)}+C_2^{(1)}-C_2^{(2)}-C_3^{(2)},\\
b_{15}=C_2^{(2)}+C_3^{(2)}-C_2^{(3)}-C_3^{(3)},\!\!\!\!&
b_{16}=C_2^{(1)}-C_1^{(2)}-C_1^{(3)}+C_3^{(3)}-C_2^{(5)}+C_1^{(6)}-2t_1+2u_1.
\end{array}
\end{eqnarray*}
The lattice $\Omega_{\Z/2\Z\times\Z/4\Z}$ admits a unique
primitive embedding in
the lattice $\Lambda_{K3}$.\\
The discriminant of $\Omega_{\Z/2\Z\times\Z/4\Z}$ is $2^{10}$ and its discriminant group is $(\Z/2\Z)^{ 2}\oplus(\Z/4\Z)^{ 4}$.\\
The group of isometries $G_{2,4}$ acts on the discriminant group
$\Omega_{\Z/2\Z\times\Z/4\Z}^{\vee}/\Omega_{\Z/2\Z\times\Z/4\Z}$ as the identity.\\
The lattice $H^2(X,\Z)^{G_{2,4}}$ is \vspace{-0.5cm}

$$
\left[\begin{array}{rrrrrr}
4& -2& 0& 0& 0& 0\\  -2& 0& -2& 0& 0& 0\\
0& -2& -64& -4& 0& 0\\  0& 0& -4& 0& -4& 0\\
0& 0& 0& -4& 80& 4\\  0& 0& 0& 0& 4& 0
\end{array}\right]
$$
and it is an overlattice of $NS(X)^{G_{2,4}}\oplus T_{X}\simeq
U(4)\oplus\left[\begin{array}{rr}-4&0\\0&-4\end{array}\right]\oplus
\left[\begin{array}{ll}4&0\\0&4\end{array}\right]$ of index two.\\

{\bf 2)} Let $G_{2,2}$ be the group generated by
$\sigma_{t_2}(=\sigma^2_{t_1})$ and
$\sigma_{u_2}(=\sigma^2_{u_1})$.\\
The invariant sublattice of the N\'eron--Severi group with respect
to $G_{2,2}$ is isometric to \label{M1}
$$\left[\begin{array}{rrrrrrrr}-4&
-4& -2& -2& -4& 0& 0& 0\\-4& -4& 0& 2& -2& 0&
  8& 4\\-2& 0& -4& 0& 0& 0& 4& 0\\-2& 2& 0& 0& 0& 0& 0& 0\\-4& -2& 0& 0& -4& 0& 6&
  0\\0& 0& 0& 0& 0& -4& 2& 0\\
0& 8& 4& 0& 6& 2& -16& 2\\0& 4& 0& 0& 0& 0& 2&
-4\end{array}\right].$$
%su mathematica ho la base... ma non so se serve.
Its orthogonal complement $(NS(X)^{G_{2,2}})^{\perp}$ is
$\Omega_{(\Z/2\Z)^2}:=(H^2(X,\mathbb{Z})^{G_{2,2}})^{\perp}$. It
is the negative definite lattice $\{\mathbb{Z}^{12}, Q\}$ where
$Q$ is the bilinear form obtained as the intersection form of the
classes
$$\begin{array}{l}
b_1=-C^{(1)}_1 - C^{(2)}_1 - C^{(3)}_1 + C^{(3)}_3 + C^{(5)}_1
+C^{(6)}_1 - 2t_1 + 2u_1,\ \ \ b_2= -C^{(1)}_1 - 2C^{(1)}_2 -
C^{(1)}_3 + F,\\ b_3= C^{(2)}_1 -
C^{(2)}_3,\ \ \ b_4= -C^{(1)}_1 - C^{(1)}_2 + C^{(2)}_1 + C^{(2)}_2,\\
b_5= C^{(1)}_1 + 2C^{(1)}_2 + C^{(1)}_3 + C^{(3)}_1 + C^{(3)}_2 +
C^{(3)}_3 + C^{(4)}_1 + C^{(4)}_2 + C^{(4)}_3 + C^{(5)}_1 +
C^{(6)}_1 -
  3F - 2s + 2u_1,\\ b_6= -C^{(1)}_1 - C^{(1)}_2 - C^{(3)}_1 - C^{(3)}_2 + F,\ \ \ b_7=C^{(4)}_1 -
  C^{(4)}_3,\ \ \
b_8= -C^{(1)}_1 - C^{(1)}_2 + C^{(4)}_2 + C^{(4)}_3,\\
b_9=C^{(5)}_1 - C^{(5)}_3,\ \ \ b_{10}=-C^{(1)}_1 - C^{(1)}_2 +
C^{(5)}_1 + C^{(5)}_2,\\ b_{11}= C^{(6)}_1 - C^{(6)}_3,\ \ \
b_{12}=
 -C^{(1)}_1 - C^{(1)}_2 + C^{(6)}_1 + C^{(6)}_2.
\end{array}
$$
The lattice $\Omega_{(\Z/2\Z)^2}$ admits a unique primitive
embedding in
the lattice $\Lambda_{K3}$.\\
The discriminant of $\Omega_{(\Z/2\Z)^2}$ is $2^{10}$ and its discriminant group is $(\Z/2\Z)^{ 6}\oplus(\Z/4\Z)^{ 2}$.\\
The group of isometries $G_{2,2}$ acts on the discriminant group
$\Omega_{(\Z/2\Z)^2}^{\vee}/\Omega_{(\Z/2\Z)^2}$ as the identity.\\
The lattice $H^2(X,\Z)^{G_{2,2}}$ is
$$
\left[\begin{array}{rrrrrrrrrr}0&4&2&0&0&2&0&-2&0&0\\
  4&0&6&-8&8&4&6&-20&8&2\\
  2&6&0&-1&2&1&2&-6&0&2\\
  0&-8&-1&-2&0&-1&2&1&2&0\\
  0&8&2&0&-4&4&0&2&0&0\\
  2&4&1&-1&4&-4&0&-1&0&-2\\
  0&6&2&2&0&0&-4&2&0&0\\
  -2&-20&-6&1&2&-1&2&-4&0&4\\
  0&8&0&2&0&0&0&0&-4&0\\
  0&2&2&0&0&-2&0&4&0&-4\end{array}\right]
$$
and it is an overlattice of $NS(X)^{G_{2,2}}\oplus T_{X}$ of index four.\\

{\bf 3)} Let $G_{4}$ be the group generated by $\sigma_{t_1}$. The
invariant sublattice of the N\'eron--Severi group with respect to
$G_{4}$ is isometric to $\langle -4\rangle\oplus
\left[\begin{array}{rrrrr}-2& 1& 0& 0& 0\\1&
-2& 4& 0& 0\\
  0& 4& 4& 8& 4\\0& 0& 8& 0& 4\\0&
   0& 4& 4& 0\end{array}\right].$
%su mathematica ho la base... ma non so se serve.
Its orthogonal complement $(NS(X)^{G_4})^{\perp}$ is
$\Omega_{\Z/4\Z}:=(H^2(X,\mathbb{Z})^{G_4})^{\perp}$. It is the
negative definite lattice $\{\mathbb{Z}^{14}, Q\}$ where $Q$ is
the bilinear form obtained as the intersection form of the classes
\begin{eqnarray*}
\begin{array}{ll}
b_1=s-t_1,&b_{i+1}=C_1^{(i)}-C_1^{(i+1)},\ i=1,2,3, \\
b_{j+3}=C_1^{(j)}-C_1^{(j)},\ j=2,\ldots, 5,&
b_{h+7}=C_2^{(h)}-C_2^{(h+1)},\ h=2,3,\\
b_{11}=C_2^{(2)}-C_1^{(1)},&b_{12}=C_2^{(1)}-C_1^{(2)},\\
b_{13}=F-C_1^{(2)}-C_2^{(2)}-C_1^{(5)}-C_2^{(5)},&
b_{14}=C_1^{(2)}+C_2^{(2)}-C_1^{(5)}-C_2^{(5)}.
\end{array}
\end{eqnarray*}
The lattice $\Omega_{\Z/4\Z}$ admits a unique primitive embedding
in
the lattice $\Lambda_{K3}$.\\
The discriminant of $\Omega_{\Z/4\Z}$ is $2^{10}$ and its discriminant group is $(\Z/2\Z)^{ 2}\oplus(\Z/4\Z)^{ 4}$.\\
The group of isometries $G_{4}$ acts on the discriminant group
$\Omega_{\Z/4\Z}^{\vee}/\Omega_{\Z/4\Z}$ as the identity.\\
The lattice $H^2(X,\Z)^{G_4}$ is
$$
\left[\begin{array}{rrrrrrrr} 0& 4& 0& 2& 0& -1& 0& 0\\  4& 0& 4&
4& -4& 0& 0&
  -4\\  0& 4 &\  0& 0& 0& 0& 0& 0\\
  2& 4& 0& 0& 0& -1& 0& 0\\  0& -4& 0& 0& -2& -1& 0&
  -2\\  -1& 0& 0& -1& -1& -2& 1& 1\\
  0& 0& 0& 0& 0& 1& -2& 0\\  0& -4& 0& 0& -2& 1& 0&
  -2
\end{array}
\right]$$ and it is an overlattice of $NS(X)^{G_{4}}\oplus T_{X}$
of index two.
\end{prop}
\bprf The proof is similar to the proof of Proposition \ref{omega and invariant 44}.\eprf

\begin{rem}{\rm The automorphisms $\sigma_{t_1}$ and $\sigma_{u_1}$ do not fix the other four elliptic fibrations described in
\cite{Keum Kondo: automorphisms of product elliptic curves}
(different from $E_1+E_2+L_1+L_2$). The involutions
$\sigma_{t_1}^2$ and $\sigma_{u_1}^2$ fix the class of the
fiber of those elliptic fibrations, however they are not induced
by torsion sections on those fibrations, in fact they do not fix
each fiber of the fibration. On the fibrations different from
$E_1+E_2+L_1+L_2$, the actions of $\sigma_{t_1}^2$ and
$\sigma_{u_1}^2$ are analogue to the ones of
$\sigma_{\mathbb{P}^1,a}$ and $\sigma_{\mathbb{P}^1,b}$ on the
fibration defined by $E_1+E_2+L_1+L_2$.}\end{rem}
\begin{prop}\label{omega and invariant 222,2222}{\bf 1)} Let $G_{2,2,2}$ be the group
generated by $\sigma_{t_2}$, $\sigma_{u_2}$ and
$\sigma_{\mathbb{P}^1,a}$. The automorphism
$\sigma_{\mathbb{P}^1,a}$ acts in the following way:
$$\begin{array}{cccc} t_1\leftrightarrow v,\ &u_1\leftrightarrow w,&
C_i^{(1)}\leftrightarrow C_i^{(2)},\ & C_i^{(5)}\leftrightarrow
C_i^{(6)},\ i=0,1,2,3,\ C_1^{(j)}\leftrightarrow C_3^{(j)},\
j=3,4,\end{array}$$ (where $w$ and $v$ are respectively the section
 $t_1+u_1+u_1$ and $t_1+t_1+u_1$ with respect to the
group law of the Mordell--Weil group) and fixes the classes $F$,
$s$ and $C_i^{(j)}$, $i=0,2$, $j=3,4$. The invariant sublattice of
the N\'eron--Severi group with respect to $G_{2,2,2}$ is isometric
to \label{M2}
$$\left[\begin{array}{rrrrrr}-4& 2& 0& 0& 0& 0\\2& -20& 6& 0& 0& 0\\0& 6& -4& -2& 0& 0\\0& 0&
-2&
    0& -2& 0\\0& 0& 0& -2& 0& -4\\0& 0& 0& 0& -4& -8\end{array}\right].$$
%su mathematica ho la base... ma non so se serve.
Its orthogonal complement $(NS(X)^{G_{2,2,2}})^{\perp}$ is
$\Omega_{(\Z/2\Z)^3}:=(H^2(X,\mathbb{Z})^{G_{2,2,2}})^{\perp}$. It
is the negative definite lattice $\{\mathbb{Z}^{14}, Q\}$ where
$Q$ is the bilinear form obtained as the intersection form of the
classes
$$
\begin{array}{l}
b_1=C^{(3)}_3-C^{(3)}_1,\ \ b_{2}=C^{(4)}_3 -C^{(4)}_1,\ \ b_{3}=
C^{(2)}_3-C^{(1)}_1,\ \ b_{4}= C^{(2)}_2-C^{(1)}_2,\ \ b_{5}=
C^{(2)}_1-C^{(1)}_1,\\
b_{6}=C^{(3)}_0 + C^{(3)}_3 -C^{(1)}_1 - C^{(1)}_2,\ \
b_{7}=C^{(3)}_1 + C^{(3)}_2 -C^{(1)}_1 - C^{(1)}_2,\ \
b_{8}=C^{(4)}_1 + C^{(4)}_2-C^{(3)}_1 - C^{(3)}_2,\\
b_{9}= C^{(5)}_2 +C^{(5)}_3 -C^{(4)}_1 - C^{(4)}_2,\ \ b_{10}=
C^{(6)}_1 + C^{(6)}_2 -C^{(5)}_2 - C^{(5)}_3,\ \
b_{11}= -C^{(1)}_1 + C^{(5)}_1 - t_1 + u_1,\\ b_{12}=C^{(6)}_1  -C^{(1)}_1 - t_1 + u_1,\ \ \ \ \ \ b_{13}= C^{(5)}_2 -C^{(1)}_2 + t_1 - u_1,\\
b_{14}=
 C^{(1)}_1 + C^{(2)}_1 + C^{(3)}_1 + C^{(4)}_1 + C^{(5)}_1 + 2C^{(5)}_2 + C^{(5)}_3 - 2F -
  2s + 2t_1.
\end{array}$$
The lattice $\Omega_{(\Z/2\Z)^3}$ admits a unique primitive
embedding in
the lattice $\Lambda_{K3}$.\\
The discriminant of $\Omega_{(\Z/2\Z)^3}$ is $2^{10}$ and its discriminant group is $(\Z/2\Z)^{ 6}\oplus (\Z/4\Z)^{ 2}$.\\
The group of isometries $G_{2,2,2}$ acts on the discriminant group
$\Omega_{(\Z/2\Z)^3}^{\vee}/\Omega_{(\Z/2\Z)^3}$ as the identity.\\
The lattice $H^2(X,\Z)^{G_{2,2,2}}$ is an overlattice of
$NS(X)^{G_{2,2,2}}\oplus T_{X}$ of index two.\\

{\bf 2)} Let $G_{2,2,2,2}$ be the group generated by
$\sigma_{t_2}$, $\sigma_{u_2}$, $\sigma_{\mathbb{P}^1,a}$ and
$\sigma_{\mathbb{P}^1,b}$. The automorphism
$\sigma_{\mathbb{P}^1,b}$ acts in the following way:
$$\begin{array}{cccc} t_1\leftrightarrow z,\ &u_1\leftrightarrow w,&
C_i^{(3)}\leftrightarrow C_i^{(4)},\ & C_i^{(5)}\leftrightarrow
C_{4-i}^{(6)},\ i=0,1,2,3,\ \ C_{1}^{(j)}\leftrightarrow
C_{3}^{(j)},\ j=1,2,
\end{array}$$ where $w$ is as in 1) and $z$ is the section
$t_1+t_1+t_1+u_1+u_1$ with respect to the group law of the
Mordell--Weil group, and fixes the classes $F$, $s$ and
$C_i^{(j)}$, $i=0,2$, $j=1,2$.  The invariant sublattice of the
N\'eron--Severi group with respect to $G_{2,2,2,2}$ is isometric
to \label{M3}
$$\left[\begin{array}{rrrrrr} -20& -8& -12& -2& 4\\-8& 8& 2& 2&
4\\-12& 2& -4& 0& 4\\-2& 2& 0& 0&
    0\\4& 4& 4& 0& -8
\end{array}\right].$$
%su mathematica ho la base... ma non so se serve.
Its orthogonal complement $(NS(X)^{G_{2,2,2,2}})^{\perp}$ is
$\Omega_{(\Z/2\Z)^4}:=(H^2(X,\mathbb{Z})^{G_{2,2,2,2}})^{\perp}$.
It is the negative definite lattice $\{\mathbb{Z}^{15}, Q\}$ where
$Q$ is the bilinear form obtained as the intersection form of the
classes
$$
\begin{array}{l}
 b_1= C^{(3)}_3-C^{(3)}_1 ,\ \ \ b_{2}=  C^{(2)}_3-C^{(1)}_1 ,\ \ \ b_{3}=
C^{(2)}_2-C^{(1)}_2 ,\ \ \ b_{4}= C^{(2)}_1-C^{(1)}_1 ,\ \ \
b_{5}= C^{(4)}_3 -C^{(3)}_1 ,\\ b_{6}= C^{(4)}_1-C^{(3)}_1  ,\ \ \
b_{7}= C^{(3)}_0 + C^{(3)}_3 -C^{(1)}_1 - C^{(1)}_2 ,\ \ \ \ \ \ \
\ \ \  b_{8}=
 C^{(3)}_1 + C^{(3)}_2-C^{(1)}_1 - C^{(1)}_2  ,\\ b_{9}= C^{(4)}_1 + C^{(4)}_2 -C^{(3)}_1 - C^{(3)}_2 ,\ \ \ b_{10}=
C^{(5)}_2 + C^{(5)}_3 -C^{(4)}_1 - C^{(4)}_2 ,\ \ \ b_{11}=
C^{(6)}_1 + C^{(6)}_2-C^{(5)}_2 - C^{(5)}_3 ,\\ b_{12}= C^{(6)}_1
-C^{(1)}_1 - t_1 + u_1,\ \ \ \ \ b_{13}= C^{(5)}_2-C^{(1)}_2  +
t_1 - u_1,\ \ \ \ \ \ \ \ \ \ b_{14}= C^{(5)}_1 -C^{(1)}_1  - t_1 + u_1,\\
b_{15}= -C^{(1)}_1 - C^{(3)}_1 - C^{(3)}_2 - C^{(4)}_1 + F +
  s - t_1.
\end{array}
$$
The lattice $\Omega_{(\Z/2\Z)^4}$ admits a unique primitive
embedding in
the lattice $\Lambda_{K3}$.\\
The discriminant of $\Omega_{(\Z/2\Z)^4}$ is $-2^{9}$ and its discriminant group is $(\Z/2\Z)^{ 6}\oplus(\Z/8\Z)$.\\
The group of isometries $G_{2,2,2,2}$ acts on the discriminant
group
$\Omega_{(\Z/2\Z)^4}^{\vee}/\Omega_{(\Z/2\Z)^4}$ as the identity.\\
The lattice $H^2(X,\Z)^{G_{2,2,2,2}}$ is an overlattice of
$NS(X)^{G_{2,2,2,2}}\oplus T_{X}$ of index two.
\end{prop}
\bprf The proof is similar  to the proof of Proposition \ref{omega
and invariant 44}, we describe here the action of
$\sigma_{\mathbb{P}^1,a}$ on the fibration (the action of
$\sigma_{\mathbb{P}^1,b}$ can be deduced in a similar way). The
automorphism $\sigma_{\mathbb{P}^1,a}$ fixes the sections $t_2$
and $u_2$, the class of the fiber (because sends fibers to fibers)
and fixes the third and the fourth reducible fiber. It acts on the
basis of the fibration $\mathbb{P}^1$ switching the two points
corresponding to the first and the second reducible fiber and the
two points corresponding to the fifth and the sixth reducible
fiber. Since the sections $t_2$ and $u_2$ are fixed, the component
$C_2^{(1)}$ (which meet the section $t_2$) is sent to the
component $C_2^{(2)}$ (which is the component of the second fiber,
which meets the same section). Similarly we obtain
$\sigma_{\mathbb{P}^1,a}(C_0^{(1)})=C_0^{(2)}$ and
$\sigma_{\mathbb{P}^1,a}(C_i^{(5)})=C_i^{(6)}$ for $i=0,2$. The
component $C_2^{(3)}$ is sent to a component of the same fiber and
since $1=C_2^{(3)}\cdot
u_2=\sigma_{\mathbb{P}^1,a}(C_2^{(3)})\cdot
\sigma_{\mathbb{P}^1,a}(u_2)=\sigma_{\mathbb{P}^1,a}(C_2^{(3)})\cdot
u_2$, we have $\sigma_{\mathbb{P}^1,a}(C_2^{(3)})=C_2^{(3)}$.
Analogously we obtain
$\sigma_{\mathbb{P}^1,a}(C_i^{(j)})=C_i^{(j)}$, $i=0,2$, $j=3,4$.
The image of $C_1^{(3)}$ has to be a component of the first fiber,
so it can only be $C_1^{(3)}$ or $C_3^{(3)}$. The curve
$C_0^{(3)}$ is a rational curve fixed by an involution, so the
involution has two fixed points on it (it cannot be fixed
pointwise because $\sigma_{\mathbb{P}^1,a}$ is a symplectic
involution and the fixed locus of a symplectic involution consists
of isolated points). These points are the intersections between
the curve $C_0^{(3)}$ and the section $u_2+t_2$ of order two and
the 0-section $s$ (which are fixed by the involution). Then the
point of intersection between $C_0^{(3)}$ and $C_1^{(3)}$ is not
fixed, and then $C_1^{(3)}$ is not fixed by the involution, so we
conclude that $\sigma_{\mathbb{P}^1,a}(C_1^{(i)})=C_3^{(i)}$,
$i=3,4$. To determine the image of the curves $C_{i}^{(j)}$,
$i=1,3$, $j=3,4,5,6$ of $u_1$ and of $t_1$ one sees that, after an
analysis of all possible combinations, the only possibility is the
action
given in the statement.%The image of the classes
%of the sections of order four can be deduced by the image of the
%classes $F$, $s$ and $C_i^{(j)}$, because they are obtained as
%linear combination with coefficient in $\Q$ of these classes.
\erem

\begin{rem}\label{rem: automorphisms P1 are induced by torsion sections}{\rm The automorphisms $\sigma_{\mathbb{P}^1,a}$, $\sigma_{\mathbb{P}^1,b}$ are
induced by a $2$-torsion section of the fibrations
$G_{22}+G_{11}+D_1+L_8$ and $G_{11}+G_{21}+E_{1}+L_{12}$
respectively (with the notation of \cite{Keum Kondo: automorphisms
of product elliptic curves}).}\end{rem}
\section{Other cases}\label{othercases}
For each group $G$ in the list \eqref{formula: symplectic group}
except for the group $G=(\Z/2\Z)^i$ for $i=3,4$ there exists a K3
surface $X$ with an elliptic fibration $\mathcal{E}_X$ such that
tors$(\MW(\mathcal{E}_X))=G$. For each group $G$ we will consider
the K3 surface $X_G$ with the minimal possible Picard number among
those admitting such an elliptic fibration. One can find the
trivial lattices of these K3 surfaces in Shimada's list (cf.
\cite[Table 2]{shimada}). In the Tables 1 and 2 in Appendix we
describe the trivial lattice of these elliptic fibrations, the
intersection properties of the torsion sections (which can be
deduced by the height formula) and we give the transcendental
lattices in all the cases except $\Z/4\Z$. In this case seems to
be more difficult to identify the transcendental lattice, however
our first aim is to compute the lattice $H^2(X,\Z)^{\Z/4\Z }$ and
we did it in Proposition \ref{omega and invariant 24 22 4}. In
Table 3 we describe $(NS(X_G)^G)^{\perp}\simeq \Omega_{G}$, giving
a basis for this lattice. The proof of the results is very similar
to the proof of Proposition \ref{NS and T case 4 4} and \ref{omega
and invariant 44} except for the equation of the fibrations.
Observe that in all the cases the Mordell--Weil group has only a
torsion part. In fact the number of moduli in the equation is
exactly $20-{\rank} ({\Tr})$, so there are no sections of infinite
order. The equation in the case of $G=\Z/2\Z$ is standard. We
computed already in Section \ref{specialauto} the equations for
$G=(\Z/2\Z)^i$, $i=2,3,4$. For $G=\Z/p\Z$, $p=3,5,7$ equations are
given e.g.\ in \cite{symplectic prime}. We now explain briefly how
to find
equations in the other cases.\\

$G=\Z/4\Z$. This fibration has in particular a section of order
two so has an equation of the kind $y^2=x(x^2+a(\tau)x+b(\tau))$,
$\deg(a(\tau))=4$, $\deg(b(\tau))=8$. To determine $a(\tau)$ and
$b(\tau)$ we observe that a smooth elliptic curve of the fibration
has a point $Q$ of order four if $Q+Q=P=(0,0)$ which is the point
of order two on each smooth fiber. Geometrically this means that
the tangent to the elliptic curve through $Q$ must intersect the
curve exactly in the point $P$. With a similar condition one
computes the equations for the groups $\Z/8\Z$,
$\Z/2\Z\times\Z/4\Z$,
$\Z/4\Z\times\Z/4\Z$.\\

$G=\Z/6\Z$. An elliptic fibration admits a 6-torsion section if
and only if it admits a 2-torsion section and a 3-torsion section,
so it is enough to impose that the fibration
$y^2=x(x^2+a(\tau)x+b(\tau))$ has an inflectional point on the
generic fiber (different from $(0:0:1)$) (cf. \cite[Ex 6,
p.38]{cassels}). By using this condition and the equation of a
fibration with two 2-torsion sections one computes the equation of
a fibration with $G=\Z/2\Z\times \Z/6\Z$.\\

$G=\Z/3\Z\times\Z/3\Z$. We consider the rational surface:
\begin{eqnarray}\label{equation Hesse pencil}
y^2=x^3+12(u_0^3u_1-u_1^4)x+2(u_0^6-20u_0^3u_1^3-8u_1^6),\,\,
(u_0:u_1)\in \PP^1
\end{eqnarray}
which has two 3-torsion sections and four singular fibers $I_3$
(cf. \cite{michela},\cite{beauville} for a description). We can
consider a double cover of $\mathbb{P}^1$ branched over two
generic points $(\alpha:1)$, $(\beta:1)$, which is realized by the
rational map $u\mapsto (\alpha\tau^2+\beta)/(\tau^2+1)$. The image
is again a copy of $\mathbb{P}^1$. This change of coordinates in
the equation \eqref{equation Hesse pencil} induces an elliptic
fibration over the second copy of $\mathbb{P}^1$. This elliptic
surface is a K3 surface, in fact the degree of the discriminant of
this Weierstrass equation is 24, and admits two 3-torsion
sections, which are the lifts of the torsion sections of
\eqref{equation Hesse pencil}. The Weierstrass equation is for
$\alpha, \beta\in \C$ \small
\begin{equation}\label{equation 3 torsion 3 torsion}
 y^2=x^3+12x(\tau^2+1)[(\alpha\tau^2+\beta)^3-(\tau^2+1)^3]
+2[(\alpha\tau^2+\beta)^6-20(\alpha\tau^2+\beta)^3(\tau^2+1)^3-8(\tau^2+1)^6].
\end{equation}

\normalsize

\vspace*{1.0cm}

In all the cases the lattices ($NS(X_G)^G)^{\perp}$ are isometric
to the lattices $\Omega_G$. Moreover, since we know the
transcendental lattice of the K3 surfaces considered (except for
$G=\Z/4\Z$), we can compute also $H^2(X_G,\Z)^G$ as in the proof
of Proposition \ref{omega and invariant 44}. Here we summarize the
main properties of the lattices $\Omega_G$ (the basis are given in
Table 3) and of their orthogonal complement:

\begin{prop}\label{prop: properties Omega and
OmegaG} For any K3 surface $X$ with a group $G$ of symplectic
automorphisms, the action on $H^2(X,\Z)$ decomposes as
$(H^2(X,\Z)^G)^{\perp_{H^2(X,\Z)}}\simeq \Omega_G$ and
$H^2(X,\Z)^G= (\Omega_G)^{\perp_{H^2(X,\Z)}}$. The following Table
lists the main properties of $\Omega_{G}$ and
$(\Omega_{G})^{\perp_{\Lambda_{K3}}}$.

\newpage

\renewcommand{\arraystretch}{1.3}

$$
\begin{array}{c|c|c|c|c|c}
G&\rank(\Omega_G)&\discr(\Omega_G)&\Omega_G^{\vee}/\Omega_G&\rank(\Omega_{G}^{\perp_{\Lambda_{K3}}})&\Omega_{G}^{\perp_{\Lambda_{K3}}}\\
\hline \Z/2\Z&8&2^8&(\Z/2\Z)^{8}&14&E_8(-2)\oplus
U^{\oplus 3}\\
\Z/3\Z&12&3^6&(\Z/3\Z)^{ 6}&10&U\oplus U(3)^{\oplus 2}\oplus A_2^{\oplus 2}\\
\Z/4\Z&14&2^{10}&(\Z/2\Z)^{2}\oplus(\Z/4\Z)^{ 4}&8&Q_4\\
\Z/5\Z&16&5^{4}&(\Z/5\Z)^{4}&6&U\oplus U(5)^{\oplus 2}\\
\Z/6\Z&16&6^{4}&(\Z/6\Z)^{4}&6&U\oplus U(6)^{\oplus 2}\\
\Z/7\Z&18&7^{3}&(\Z/7\Z)^{3}&4&U(7)\oplus \left[\begin{array}{ll}4&1\\1&2\end{array}\right]\\
\Z/8\Z&18&8^{3}&\Z/2\Z\oplus\Z/4\Z\oplus (\Z/8\Z)^{2}&4&U(8)\oplus \left[\begin{array}{ll}2&0\\0&4\end{array}\right]\\
(\Z/2\Z)^2&12&2^{10}&(\Z/2\Z)^{6}\oplus(\Z/4\Z)^2&10&U(2)^{\oplus 2}\oplus Q_{2,2}\\
(\Z/2\Z)^3&14&2^{10}&(\Z/2\Z)^{6}\oplus(\Z/4\Z)^2&8&U(2)^{\oplus
3}
\oplus \langle -4\rangle^{\oplus2}\\
(\Z/2\Z)^4&15&-2^{9}&(\Z/2\Z)^{6}\oplus\Z/8\Z&7&\langle -8\rangle\oplus U(2)^{\oplus 3}\\
\Z/2\Z\times\Z/4\Z&16&2^{10}&(\Z/2\Z)^{2}\oplus(\Z/4\Z)^4&6&Q_{2,4}\\
\Z/2\Z\times\Z/6\Z&18&2^{4}3^3&(\Z/3\Z)\oplus(\Z/12\Z)^2&4&\left[\begin{array}{rrrr}0& 6& 0& 0\\ 6& 0& -3& 0\\ 0& -3& 6& 6\\
0& 0& 6& 8
\end{array}\right]\\
(\Z/3\Z)^2&16&3^{6}&(\Z/3\Z)^{4}\oplus \Z/9\Z&6&U(3)^{\oplus 2}\oplus \left[\begin{array}{ll}2&3\\3&0\end{array}\right]\\
(\Z/4\Z)^2&18&2^{8}&(\Z/2\Z)^{2}\oplus (\Z/8\Z)^2&4&\left[\begin{array}{rrrr}
  4& 6& 0& 0\\  6& 4& 6& 4\\  0& 6& 4& 0\\
  0& 4& 0& 0\end{array}\right]
\end{array}
$$

\renewcommand{\arraystretch}{1.0}

where\\
$Q_4=\left[
\begin{array}{rrrrrrrr}
0& 4& 0& 2& 0& -1& 0& 0\\  4& 0& 4& 4& -4& 0& 0&
  -4\\  0& 4& 0& 0& 0& 0& 0& 0\\
  2& 4& 0& 0& 0& -1& 0& 0\\  0& -4& 0& 0& -2& -1& 0&
  -2\\  -1& 0& 0& -1& -1& -2& 1& 1\\
  0& 0& 0& 0& 0& 1& -2& 0\\  0& -4& 0& 0& -2& 1& 0&
  -2
\end{array}\right],$\\
$Q_{2,2}=\left[\begin{array}{rrrrrr}0& 1& 0& 0& 0& 0\\
1& -2& 2& 0& 0& 0\\
0& 2& -4& 2& 0& 0\\
0& 0& 2& -4& 2& 0\\
0& 0& 0& 2& -4& 4\\
0& 0& 0& 0& 4& -8 \end{array}\right],\
Q_{2,4}=\left[\begin{array}{rrrrrr}
4& -2& 0& 0& 0& 0\\  -2& 0& -2& 0& 0& 0\\
0& -2& -64& -4& 0& 0\\  0& 0& -4& 0& -4& 0\\
0& 0& 0& -4& 80& 4\\  0& 0& 0& 0& 4& 0
\end{array}\right].\\
$ The lattices $\Omega_G$ are even negative definite lattices, they do
not contain vectors of length $-2$ and are generated by vectors of
maximal length, i.e. by vectors of length $-4$.
\end{prop}
\bprf The proof is for each case similar to the proof of
Proposition \ref{omega and invariant 44}. The fact that the
lattices $\Omega_G$ do not contain vector of length $-2$ was
proved by Nikulin (cf. \cite{Nikulin symplectic}). The fact that
these lattices are generated by vectors of length $-4$ is a direct
consequence of the choice of the basis in Table 3 of the Appendix,
and in Propositions \ref{omega and invariant 44}, \ref{omega and
invariant 24 22 4}, \ref{omega and invariant 222,2222}. In fact
all these bases consist of classes with self intersection
$-4$.\erem
\begin{rem}{\rm The lattices $\Omega_{(\Z/2\Z)^2}$, $\Omega_{\Z/4\Z}$, $\Omega_{(\Z/2\Z)^4}$ and $\Omega_{(\Z/3\Z)^2}$ are isometric respectively to the
laminated lattices (with the bilinear form multiplied by $-1$)  $\Lambda_{12}(-1)$, $\Lambda_{14.3}(-1)$,
$\Lambda_{15}(-1)$ and to the lattice $K_{16.3}(-1)$ (a special sublattice of the Leech-lattice). All these lattices are described in \cite{Plesken Pohst}, \cite{Conway Sloane} and \cite{nebe} and the isometries
can be proved applying an algorithm described in
\cite{PS} about isometries of lattices. We
did the computations by using the package {\it Automorphism group
and isometry testing} of MAGMA (which uses the algorithm of \cite{PS})
and the computation can be done at the
web page MAGMA-Calculator (cf. \cite{magma}).\\
In particular the lattices $\Lambda_{12}(-1)$, $\Lambda_{14.3}(-1)$,
$\Lambda_{15}(-1)$ and $K_{16.3}(-1)$ are primitive sublattices of the lattice
$\Lambda_{K3}$.\\
The lattice $\Omega_{\Z/3\Z}$ is isometric to the lattice
$K_{12}(-2)$ described in \cite{Conway Sloane k12} and
\cite{Coxeter Todd}. The proof of this isometry can be found in
\cite{symplectic prime}, where a description of the lattices
$\Omega_{\Z/5\Z}$ and $\Omega_{\Z/7\Z}$ is also given.}
\end{rem}
\begin{rem}{\rm
We have computed here for each group of symplectic automorphisms
in the list of Nikulin the invariant sublattice and its orthogonal
complement in the K3 lattice. It is in general difficult to find
explicitly the action of a symplectic automorphism on
$\Lambda_{K3}$. This is done by Morrison \cite{morrison} in the
case of involutions, otherwise this is not known. By using an
elliptic fibration with a fiber $I_{16}$, a fiber $I_4$, four
fibers $I_1$ and a 4-torsion section $s_1$ (cf.
\cite[No.3171]{shimada}) one can identify an operation of a
symplectic automorphism of order four on the K3 lattice whose
square is the operation described by Morrison. In fact it is easy
to find two copies of the lattice $E_8(-1)$ in the fiber $I_{16}$
together with the sections $s_0$, $s_2$. Then the operation of
$s_2$ interchanges these two lattices and after the identification
with $\Lambda_{K3}$, it fixes a copy of $U^3$. Since the
computations are quite involved we omit them.}
\end{rem}
\section{Families of K3 surfaces admitting symplectic automorphisms}\label{families}
Now let us consider an algebraic K3 surface $X$ with a group of
symplectic automorphisms $G$. Since the action of $G$ is trivial
on the transcendental lattice, we have $T_X\subset H^2(X,\Z)^G$
and $NS(X)\supset (H^2(X,\Z)^G)^{\perp}\simeq \Omega_G$. The
lattice $\Omega_G$ is negative definite, and the signature of
$NS(X)$ is $(1,\rho-1)$, so $\Omega_G^{\perp_{NS(X)}}$ has to
contain a class with positive self-intersection. In particular the
minimal possible Picard number of $X$ is $\rank(\Omega_G)+1$. Here
we want to consider the possible N\'eron--Severi groups of K3
surfaces admitting a finite abelian group of symplectic
automorphisms. The first step is to find all lattices
$\mathcal{L}$ such that: $\rank(\mathcal{L})=\rank(\Omega_G)+1$,
$\Omega_G$ is primitively embedded in $\mathcal{L}$ and
$\Omega_G^{\perp_{\mathcal{L}}}$ is positive definite. We prove
that if a K3 surface admits one of these lattices as
N\'eron--Severi group, then it admits the group $G$ as symplectic
group of automorphisms and viceversa if a K3 surface with Picard
number $\rank(\Omega_G)+1$ admits $G$ as symplectic group of
automorphisms, then its N\'eron--Severi group is isometric to one
of those lattices.\\
To find the lattice $\mathcal{L}$ we need the following
remark.
\begin{rem}\label{rem: even overlattice adding a vector}
{\rm We recall here the correspondence between even overlattices
of a lattice and totally isotropic subgroups of its discriminant
group (cf. \cite[Proposition 1.4.1]{Nikulin bilinear}). Let $L_H$
be an overlattice of $L$ such that $\rank (L)=\rank(L_H)$. If
$L\hookrightarrow L_H$ with a finite index, then
$[L_H:L]^2=d(L)/d(L_H)$ and $L\subseteqq L_H\subseteqq L^{\vee}$.
The lattice $L_H$ corresponds to a subgroup of the discriminant
group $H$ of $L^{\vee}$. Viceversa a subgroup $H$ of $L^{\vee}/L$
corresponds to a $\Z$-module $L_H$ such that $L\subseteqq L_H$. On
$L^{\vee}$ there is defined a bilinear form, which is the
$\Q$-linear extension of the bilinear form on $L$ and so on every
subgroup of $L^{\vee}$ there is defined a bilinear form, induced
by the form on $L$.\\
Let $L$ be an even lattice with bilinear form $b$ and let $H$ be a
subgroup of $L^{\vee}/L$. Denote also by $b$ the form induced on
$H$ and on $L_H$ by the form on $L$.
If $H$ is such that $b(h,h)\in 2\Z$ for each $h\in H$, then $L_H$ is an even overlattice of $L$.\\
\binf each element $x\in L_H$ can be written as $x=h+l$, $h\in H$,
$l\in L$. The value of the quadratic form $q(x)=b(x,x)$ on $x$ is
in $2\Z$, in fact $b(x,x)=b(h+l,h+l)=b(h,h)+2b(h,l)+b(l,l)$ and
$b(h,h),b(l,l)\in 2\Z$ by hypothesis, $b(h,l)\in\Z$ because $h\in
L^{\vee}$. For each $h,h'\in H$ we have
$b(h+h',h+h')=b(h,h)+2b(h,h')+b(h',h')\in 2\Z$ by hypothesis and
so $b(h,h')\in\Z$. Moreover,
$b(x,x')=b(h+l,h'+l')=b(h,h')+b(h,l')+b(l,h')+b(l,l')$ and
$b(h,l')$, $b(l,h')\in\Z$ because $h\in L^{\vee}$, $b(l,l')\in
\Z$, because $l,l'\in L$, and we have already proved that
$b(h,h')\in\Z$. So the bilinear form on $L_H$ takes values in $\Z$
and the quadratic form induced by $b$ takes values in $2\Z$, i.e.
$L_H$ is an even lattice.\erem}\end{rem}

\begin{prop}\label{prop: overlattice of L+Upsilon}
Let $X$ be an algebraic K3 surface and let $\Upsilon$ be a negative
definite primitive sublattice of $NS(X)$ such that
$\Upsilon^{\perp_{NS(X)}}=\Z L$ where $L$ is a class in $NS(X)$,
$L^2=2d$ with $d>0$. Let the discriminant group of $\Upsilon$ be
$\oplus_{j=1}^m (\Z/h_j\Z)^{n_j}$, with $h_j|h_{j+1}$ for each
$j=1,\ldots, m-1$ and put $\Z L\oplus \Upsilon=\mathcal{L}$.
Then $NS(X)$ is one of the following lattices:\\
{\it i)} if $\gcd(2d,h_j)=1$ for all $j$, then
$NS(X)=\mathcal{L}$;\\
{\it ii)} if $L^2\equiv 0 \mod r$ with
$r|h_m$, then either $NS(X)=\mathcal{L}$, or
$NS(X)=\mathcal{L}'_r$ is an overlattice of $\mathcal{L}$ of index
$r$ with $\Upsilon$ primitively embedded in $\mathcal{L}'_r$ . If
$NS(X)=\mathcal{L}'_{r}$ then there is an element of type
$(L/r,\upsilon/r)\in \mathcal{L}'_{r}$ which is not in
$\mathcal{L}$, where $\upsilon/r\in\Upsilon^{\vee}$,
$\upsilon^2\equiv -2d\mod2r^2$.
\end{prop}
\bprf Let us suppose that $L^2\equiv 0\mod r$, where $r|h_i$ for a
certain $i$. If there exists an element $\upsilon$ in $\Upsilon$ such that
$\upsilon/r\in \Upsilon^{\vee}$ and $\upsilon^2\equiv -2d\mod
2r^2$, then the element $v=(L/r,\upsilon/r)$ has an integer
intersection with all the classes of $\mathcal{L}$ and has an even
self-intersection. Then (cf. Remark \ref{rem: even overlattice
adding a vector}) adding $v$ to the lattice $\mathcal{L}$ we find
an even
overlattice of $\mathcal{L}$ of index $r$.\\
We prove that if there exists an even overlattice $\mathcal{L}_r'$
of $\mathcal{L}$, with $\Upsilon$ primitively embedded in
$\mathcal{L}_r'$, then there is a $j$ with $\gcd(2d,h_j)>1$ (i.e.
there exists an integer $r$ such that $L^2\equiv 0\mod r$ and
$r|h_m$) and the overlattice $\mathcal{L}_r'$ is constructed by
adding a class $(L/r,v/r)$ to $\mathcal{L}$. Let $(\alpha
L/r,\beta\upsilon'/s)$, $\upsilon'\in \Upsilon$, $\alpha,$
$\beta$, $r$, $s\in\Z$ be an element in $NS(X)$, then its
intersection with $L$ and with the classes of $\Upsilon$ is an
integer. This implies that $\alpha L/r$ and $\beta \upsilon'/s$
are classes in $(\Z L\oplus \Upsilon)^{\vee}$, w.l.o.g we assume
that $\gcd(\alpha,r)=1$, $\gcd(\beta,s)=1$. In particular $\alpha
L/r=k L/(2d)\in (\Z L)^{\vee}$ for a certain $k\in\Z$ and
$\beta\upsilon'/s\in\Upsilon^{\vee}$ so $r|2d$ and $s|h_m$.\\
By the relations
$$
\begin{array}{l}
r(\alpha L/r,\beta\upsilon'/s)-\alpha
L=r\beta\upsilon'/s\in NS(X),\\
s(\alpha L/r,\beta\upsilon'/s)-\beta \upsilon'=s\alpha L/r\in
NS(X)
\end{array}
$$
and the fact that $\Upsilon$ is a primitive sublattice of $NS(X)$,
we obtain that $r=s$. The class is now $(\alpha
L/r,\beta\upsilon'/r)$ with $r|2d$ (so $L^2=2d\equiv 0 \mod r$)
and $r|h_m$. Since $\gcd(\alpha,r)=1$, there exist $a,b\in\Z$ such
that $a\alpha+br=1$, so the class $a(\alpha L/r,\beta
\upsilon'/r)-bL=(L/r, a\beta \upsilon'/r)$ is in the
N\'eron--Severi group. Hence we can assume $\alpha=1$ and we take
the class $(L/r,\upsilon/r)$, where $\upsilon=a\beta\upsilon'$$\in
\Upsilon$. Observe that the self-intersection of a class in
$NS(X)$ is an even integer, so
$(L/r,\upsilon/r)^2=(2d+\upsilon^2)/r^2\in 2\Z$, hence
$\upsilon^2\equiv -2d\mod 2r^2$.\\
Now we prove that all the possible overlattices of $\Z L\oplus
\Upsilon$ containing $\Upsilon$ as a primitive sublattice, are
obtained in this way. Essentially we need to prove that if there
exists an overlattice $\mathcal{L}'$  of a lattice
$\mathcal{L}'_r$ (containing $\Upsilon$ primitively), then it is a
lattice $\mathcal{L}'_{r'}$ with $r|r'$, so all the possible
overlattices $\mathcal{L}'$ of $\Z L\oplus \Upsilon$ with
$\Upsilon$ primitive in $\mathcal{L}'$ are of type
$\mathcal{L}'_k$ for a certain $k$. Assume that
$(L/t,\upsilon/t)$, $(L/s,w/s)$ and $\Z L\oplus \Upsilon$ generate
the lattice $\mathcal{L}'$. If $t=s$, then we have $(\upsilon
-w)/t$$\in NS(X)$ and so $\cl'$ is generated by $(L/t,\upsilon/t)$
and $\cl$. If $s\not=t$ consider the element
$(L/t,\upsilon/t)-(L/s,w/s)=((s-t)L/ts,(s\upsilon-wt)/(ts))$, if
gcd$(s-t,ts)=1$, then we can replace
$((s-t)L/ts,(s\upsilon-wt)/(ts))$ by an element
$(L/ts,\upsilon'/ts)$. Then it is an easy computation to see that
the lattice $\cl'$ is generated by this element and $\cl$. If
$\gcd((s-t),ts)\neq 1$, then we reduce the fraction $(t-s)/ts$ and
we apply the same arguments.\erem

\begin{prop}\label{prop: possible NS of K3 with G symplectic} Let $X_G$ be an algebraic K3 surface
with a finite abelian group $G$ as group of symplectic
automorphisms and with $\rho(X)=\rank(\Omega_G)+1$. Then
the N\'eron--Severi group of $X$ is one of the following (we write $a\equiv_d b$ for $a\equiv b\mod d$) \small\\
\begin{eqnarray*}
\begin{array}{l}
\bullet\ G=\\
\Z/2\Z:\end{array}\left\{\begin{array}{l} L^2\equiv_4 0 ,\
\left\{\begin{array}{l}NS(X)=L\oplus \Omega_G\mbox{ \rm or
}\\NS(X)/(L\oplus \Omega_G)=\langle(L/2,v/2)\rangle,\
v/2\in\Omega_G^{\vee}/\Omega_G,
L^2+v^2\equiv_8 0,\end{array}\right.\\
 L^2\not\equiv_4 0,\   NS(X)=L\oplus
\Omega_G;\end{array}\right.
\end{eqnarray*}
\begin{eqnarray*}
\begin{array}{l}\bullet\ G=\\
\Z/p\Z\\
p=3,5,7\end{array}: \left\{\begin{array}{l} L^2\equiv_{2p} 0 ,\
\left\{\begin{array}{l}NS(X)=L\oplus \Omega_G\mbox{ \rm or
}\\NS(X)/(L\oplus \Omega_G)=\langle(L/p,v/p)\rangle,\
v/p\in\Omega_G^{\vee}/\Omega_G,  L^2+v^2\equiv_{2p^2}
0,\end{array}\right.\\ L^2\not\equiv_{2p} 0,\ NS(X)=L\oplus
\Omega_G;\end{array}\right.\end{eqnarray*}
\begin{eqnarray*}
\begin{array}{l}\bullet\ G=\\
\Z/4\Z\\
(\Z/2\Z)^2\\
(\Z/2\Z)^3\\
\Z/2\Z\times\Z/4\Z\end{array}: \left\{\begin{array}{l}
L^2\equiv_{4}\! 0 ,\  \left\{\begin{array}{l}NS(X)=L\oplus
\Omega_G\mbox{ \rm or }\\NS(X)/(L\oplus
\Omega_G)=\langle(L/2,v/2)\rangle,\
v/2\in\Omega_G^{\vee}/\Omega_G,  L^2+v^2\equiv_{8}
0,\\NS(X)/(L\oplus \Omega_G)=\langle(L/4,v/4)\rangle,\
v/4\in\Omega_G^{\vee}/\Omega_G,  L^2+v^2\equiv_{32}
0,\end{array}\right.\\ L^2\not\equiv\!_{4} 0,\
\left\{\begin{array}{l}NS(X)=L\oplus \Omega_G\mbox{ \rm or
}\\NS(X)/(L\oplus \Omega_G)=\langle(L/2,v/2)\rangle,\
v/2\in\Omega_G^{\vee}/\Omega_G,  L^2+v^2\equiv_8
0;\end{array}\right.\end{array}\right.\end{eqnarray*}
\begin{eqnarray*}
\begin{array}{l}\bullet\ G=\\
\Z/6\Z\end{array}: \left\{\begin{array}{l} L^2\equiv_{6} 0 ,\
\left\{\begin{array}{l}NS(X)=L\oplus \Omega_G\mbox{ \rm or
}\\NS(X)/(L\oplus \Omega_G)=\langle(L/2,v/2)\rangle,\
v/2\in\Omega_G^{\vee}/\Omega_G,  L^2+v^2\equiv_{8}
0,\\NS(X)/(L\oplus \Omega_G)=\langle(L/3,v/3)\rangle,\
v/3\in\Omega_G^{\vee}/\Omega_G,  L^2+v^2\equiv_{18} 0,
\\NS(X)/(L\oplus \Omega_G)=(L/6,v/6), \ v/6\in\Omega_G^{\vee}/\Omega_G,  \
L^2+v^2\equiv_{72} 0,
\end{array}\right.\\ L^2\not\equiv\!_{6}
0,\   \left\{\begin{array}{l}NS(X)=L\oplus \Omega_G\mbox{ \rm or
}\\NS(X)/(L\oplus \Omega_G)=(L/2,v/2),\
v/2\in\Omega_G^{\vee}/\Omega_G,  L^2+v^2\equiv_8
0;\end{array}\right.\end{array}\right.\end{eqnarray*}
\begin{eqnarray*}
\begin{array}{l}\bullet\ G=\\
\Z/8\Z\\
(\Z/2\Z)^4\\
(\Z/4\Z)^2\end{array}: \left\{\begin{array}{l} L^2\equiv_{8} 0 ,\
\left\{\begin{array}{l}NS(X)=L\oplus \Omega_G\mbox{  \rm or
}\\NS(X)/(L\oplus \Omega_G)=\langle(L/2,v/2)\rangle,\
v/2\in\Omega_G^{\vee}/\Omega_G,  L^2+v^2\equiv_{8}
0,\\NS(X)/(L\oplus \Omega_G)=\langle(L/4,v/4)\rangle,\
v/4\in\Omega_G^{\vee}/\Omega_G,  L^2+v^2\equiv_{32} 0,
\\NS(X)/(L\oplus \Omega_G)=(L/8,v/8), \ v/8\in\Omega_G^{\vee}/\Omega_G,  \
L^2+v^2\equiv_{128} 0,
\end{array}\right.\\
L^2\equiv_{4} 0,\  L^2\not\equiv_8 0 ,\
\left\{\begin{array}{l}NS(X)=L\oplus \Omega_G\mbox{ \rm  or
}\\NS(X)/(L\oplus \Omega_G)=\langle(L/2,v/2)\rangle,\
v/2\in\Omega_G^{\vee}/\Omega_G,  L^2+v^2\equiv_{8}
0,\\NS(X)/(L\oplus \Omega_G)=\langle(L/4,v/4)\rangle,\
v/4\in\Omega_G^{\vee}/\Omega_G,  L^2+v^2\equiv_{32} 0,
\end{array}\right.\\
L^2\not\equiv\!_{4} 0,\ \left\{\begin{array}{l}NS(X)=L\oplus
\Omega_G\mbox{\rm or }\\NS(X)/(L\oplus \Omega_G)=(L/2,v/2),\
v/2\in\Omega_G^{\vee}/\Omega_G,  L^2+v^2\equiv_8
0;\end{array}\right.\end{array}\right.\end{eqnarray*}
\begin{eqnarray*}
\begin{array}{l}\bullet\ G=\\
(\Z/3\Z)^2\\ \end{array}: \left\{\begin{array}{l} L^2\equiv_{18} 0
,\  \left\{\begin{array}{l}NS(X)=L\oplus \Omega_G\mbox{ \rm or
}\\NS(X)/(L\oplus \Omega_G)=\langle(L/3,v/3)\rangle,\
v/3\in\Omega_G^{\vee}/\Omega_G,  L^2+v^2\equiv_{18}
0,\\NS(X)/(L\oplus \Omega_G)=\langle(L/9,v/9)\rangle,\
v/9\in\Omega_G^{\vee}/\Omega_G,  L^2+v^2\equiv_{162} 0,
\end{array}\right.\\ L^2\not\equiv\!_{18}0,\ L^2\equiv_{6} 0
,\   \left\{\begin{array}{l}NS(X)=L\oplus \Omega_G\mbox{\rm or
}\\NS(X)/(L\oplus \Omega_G)=(L/3,v/3),\
v/3\in\Omega_G^{\vee}/\Omega_G,  L^2+v^2\equiv_{18}
0,\end{array}\right.\\
L^2\not\equiv\!_{6}0,\ NS(X)=L\oplus \Omega_G;
\end{array}\right.\end{eqnarray*}
\begin{eqnarray*}
\begin{array}{l}\bullet\ G=\\
\Z/6\Z\times\Z/2\Z\end{array}: \left\{\begin{array}{l}
L^2\equiv_{12} 0 ,\  \left\{\begin{array}{l}NS(X)=L\oplus
\Omega_G\mbox{ \rm or }\\NS(X)/(L\oplus
\Omega_G)=\langle(L/12,v/12)\rangle,\
v/12\in\Omega_G^{\vee}/\Omega_G,  L^2+v^2\equiv_{288}
0,\\NS(X)/(L\oplus \Omega_G)=\langle(L/6,v/6)\rangle,\
v/6\in\Omega_G^{\vee}/\Omega_G,  L^2+v^2\equiv_{72} 0,\\
NS(X)/(L\oplus \Omega_G)=\langle(L/4,v/4)\rangle,\
v/4\in\Omega_G^{\vee}/\Omega_G,  L^2+v^2\equiv_{32} 0,\\
NS(X)/(L\oplus \Omega_G)=\langle(L/3,v/3)\rangle,\
v/3\in\Omega_G^{\vee}/\Omega_G,  L^2+v^2\equiv_{18} 0,\\
NS(X)/(L\oplus \Omega_G)=\langle(L/2,v/2)\rangle,\
v/2\in\Omega_G^{\vee}/\Omega_G,  L^2+v^2\equiv_{8} 0,
\end{array}\right.\\
 L^2\not\equiv_{12} 0, L^2\equiv_6 0  ,\
\left\{\begin{array}{l}NS(X)=L\oplus \Omega_G\mbox{ \rm or
}\\NS(X)/(L\oplus \Omega_G)=\langle(L/6,v/6)\rangle,\
v/6\in\Omega_G^{\vee}/\Omega_G,  L^2+v^2\equiv_{72} 0,\\
NS(X)/(L\oplus \Omega_G)=\langle(L/3,v/3)\rangle,\
v/3\in\Omega_G^{\vee}/\Omega_G,  L^2+v^2\equiv_{18} 0,\\
NS(X)/(L\oplus \Omega_G)=\langle(L/2,v/2)\rangle,\
v/2\in\Omega_G^{\vee}/\Omega_G,  L^2+v^2\equiv_{8} 0,
\end{array}\right.\\
 L^2\not\equiv_{12} 0, L^2\equiv_4 0  ,\
\left\{\begin{array}{l}NS(X)=L\oplus \Omega_G\mbox{\rm or
}\\NS(X)/(L\oplus \Omega_G)=\langle(L/4,v/4)\rangle,\
v/4\in\Omega_G^{\vee}/\Omega_G,  L^2+v^2\equiv_{32} 0,\\
NS(X)/(L\oplus \Omega_G)=\langle(L/2,v/2)\rangle,\
v/2\in\Omega_G^{\vee}/\Omega_G,  L^2+v^2\equiv_{8} 0,
\end{array}\right.\\
 L^2\not\equiv\!_{6}0,\ L^2\not\equiv_{4} 0 ,\
\left\{\begin{array}{l}NS(X)=L\oplus \Omega_G\mbox{ \rm or
}\\NS(X)/(L\oplus \Omega_G)=(L/2,v/2),\
v/2\in\Omega_G^{\vee}/\Omega_G,  L^2+v^2\equiv_{8}
0.\end{array}\right.\end{array}\right.\end{eqnarray*} \normalsize
Moreover the class $L$ can be chosen as an ample class. We will
denote by $\mathcal{L}_G^{2d}$ the lattice $\Z L\oplus \Omega_G$,
where $L^2=2d$, and by $\mathcal{L}_{G,r}^{'2d}$ the
overlattices of $\mathcal{L}_G^{2d}$ of index $r$.\\
\end{prop}
\bprf The statement follows from Proposition
\ref{prop: overlattice of L+Upsilon} for to the lattices $\Z L\oplus
\Omega_G$, where the discriminant group of $\Omega_G$ is given in
Proposition \ref{prop: properties Omega and OmegaG}. The cases
$G=\Z/p\Z$, $p=3,5,7$ are described in \cite{symplectic prime}.
The case $G=\Z/2\Z$ is described in \cite{bert Nikulin
involutions}. Applying Proposition \ref{prop: overlattice of
L+Upsilon} one obtains that there are no conditions on $L^2$ to
obtain an overlattice of index two of $\Z L\oplus \Omega_G$.
However if $G=\Z/2\Z$ and $L^2\not\equiv_4 0$ there are no
possible overlattices of $\Z L\oplus \Omega_{\Z/2\Z}$. This
depends on the properties of $\Omega_{\Z/2\Z}$. In fact if there
were an overlattice of $\Z L\oplus \Omega_{\Z/2\Z}$,
$L^2\not\equiv_4 0 $, then there would be an element of the form
$(L/2,v/2)$ such that $L^2+v^2\equiv 0\mod 8$. Since
$L^2\not\equiv_4 0$ then $L^2\equiv_8 2$ or $L^2\equiv_8 6$. This
implies that $v^2\equiv_8 \pm 2$. But all the elements in
$\Omega_{\Z/2\Z}=E_8(-2)$ have self-intersection which is a
multiple of four.\\
Since $L^2>0$, by the Riemann--Roch theorem we can assume that $L$
or $-L$ effective. Hence we assume $L$ to be effective. Let $D$ be
an effective $(-2)$-curve, then $D=\alpha L+v'$, with
$v'\in\Omega_{G}$ and $\alpha>0$ since $\Omega_{G}$ does not
contain $(-2)$-curves. We have $L\cdot D=\alpha L^2>0$, and so $L$
is ample. \erem

\begin{prop}\label{prop: if LG is Ns then G acts symplectically}
Let $\mathcal{L}_G$ be either $\mathcal{L}_G^{2d}$, or
$\mathcal{L}_{G,r}^{'2d}$. Let $X_G$ be an algebraic K3 surface
such that $NS(X_G)=\mathcal{L}_G$. Then $X_G$ admits $G$ as group
of symplectic automorphisms (cf. also \cite[Theorem 4.15]{Nikulin
symplectic}, \cite[Proposition 2.3]{bert Nikulin involutions},
\cite[Proposition 5.2]{symplectic prime}).
\end{prop}
\bprf
Denote by $\widetilde{G}$ a group of isometries of $E_8(-1)^2\oplus U^3$ acting as in the
case of the K3-surfaces with an elliptic fibration and a group of symplectic
automorphisms $G$.\\
\textit{Step 1: the isometries of $\widetilde{G}$ fix the
sublattice $\mathcal{L}_G$.} Since
$\widetilde{G}(\Omega_G)=\Omega_G$ and $\widetilde{G}(L)=L$
(because $L\in\Omega_G^{\perp}$ which is the invariant sublattice
of $\Lambda_{K3}$), if $\mathcal{L}_G=\mathcal{L}^{2d}_{G}=\Z
L\oplus \Omega_G$ it is clear that
$\widetilde{G}(\mathcal{L}_G)=\mathcal{L}_G$. Now we consider the
case $\mathcal{L}_G=\cl'^{2d}_{G}$. The isometry $\widetilde{G}$
acts trivially on $\Omega_G^{\vee}/\Omega_G$ (it can be proved by
a computer-computation on the generators of the discriminant form
of the lattices $\Omega_G$ and on
$(\mathbb{Z}L)^{\vee}/\mathbb{Z}L$).\\ Let $\frac{1}{r}(L,v')\in
\mathcal{L}_G$, with $v'\in\Omega_G$. This is also an element in
$(\Omega_G\oplus L\mathbb{Z})^{\vee}/(\Omega_G\oplus
L\mathbb{Z})$. So we have
$\widetilde{G}(\frac{1}{r}(L,v'))\equiv\frac{1}{r}(L,v')\mod(\Omega_G\oplus\mathbb{Z}L)$,
which means
$\widetilde{G}(\frac{1}{r}(L,v'))=\frac{1}{r}(L,v')+(\beta
L,v'')$, $\beta \in\mathbb{Z},$ $v''\in\Omega_G.$
Hence in any case we have $\widetilde{G}(\mathcal{L}_G)=\mathcal{L}_G$.\\
\textit{Step 2: The isometries of $\widetilde{G}$ preserve the
K\"ahler cone $\mathcal{C}_X^+$.} We recall that the K\"ahler cone
of a K3 surface $X$ can be described as the set
$$\mathcal{C}_X^+=\{x\in V(X)^+\ |\ (x,d)>0\ \textrm{for each}\ d\in NS(X)\
\textrm{such that}\ (d,d)=-2,\ d\ \textrm{effective}\}$$ where
$V(X)^+$ is the connected component of $\{x\in
H^{1,1,}(X,\mathbb{R})\,|\, (x,x)>0\}$ containing a K\"ahler
class. First we prove that $\widetilde{G}$ fixes the set of the
effective $(-2)$-classes. Since there are no $(-2)$-classes in
$\Omega_G$, if $N\in\mathcal{L}_G$ has $N^2=-2$, then
$N=\frac{1}{r}(aL,v')$, $v'\in\Omega_G$, for an integer $a\neq 0$
(recall that $r\in\Z_{>0}$). Since $\frac{1}{r} aL^2=L\cdot N>0$,
because $L$ and $N$ are effective divisors, we obtain $a>0$. The
curve $N'=\widetilde{G}(N)$ is a $(-2)$-class because
$\widetilde{G}$ is an isometry, hence $N'$ or $-N'$ is effective.
Since $N'=\widetilde{G}(N)=(aL/r,\widetilde{G}(v')/r)$ we have
$-N'\cdot L=-a/r L^2<0$ and so $-N'$ is not effective. Hence
$N'=\widetilde{G}(N)$ is an effective $(-2)$-class. Using the fact
that $\widetilde{G}$ has finite order, it is clear that
$\widetilde{G}$ fixes the set of the effective $(-2)$-classes.\\
Now let $x\in\mathcal{C}_X^+$, then
$\widetilde{G}(x)\in\mathcal{C}_X^+$, in fact
$(\widetilde{G}(x),\widetilde{G}(x))=(x,x)>0$ and for each
effective $(-2)$-class $d$ there exists an effective
$(-2)$-class $d'$ with $d=\widetilde{G}(d')$, so we have
$(\widetilde{G}(x),d)=(\widetilde{G}(x),\widetilde{G}(d'))=(x,d')>0$.
Hence $\widetilde{G}$ preserves $\mathcal{C}_X^+$ as claimed.\\
\textit{Step 3: The isometries of $\widetilde{G}$ are induced by
automorphisms of the surface $X$.} The isometries of
$\widetilde{G}$ are the identity on the sublattice
$\mathcal{L}_G^{\perp}$ of $\Lambda_{K3}$, so they are Hodge
isometries. By the Torelli theorem, an effective Hodge isometry of
the lattice $\Lambda_{K3}$ is induced by an automorphism of the K3
surface (cf. \cite[Theorem 1.11]{bpv}). By
\cite[Corollary 3.11]{bpv}, $\widetilde{G}$ preserves the set of
the effective divisors if and only if it preserves the K\"ahler
cone. The isometries of $\widetilde{G}$ preserve the K\"ahler
cone, from Step 2, and so they are effective, hence induced by automorphisms of
the surface. The latter are symplectic by construction (they are the
identity on the transcendental lattice
$T_X\subset\Omega_G^{\perp}$).\eprf

The previous proposition implies the following:
\begin{prop} Let $\mathcal{L}_G$ be either $\mathcal{L}_G^{2d}$, or
$\mathcal{L}_{G,r}^{'2d}$. Let $X$ be an algebraic K3 surface with
Picard number equal to $1+\rank(\Omega_G)$ for a certain $G$
in the list \eqref{formula: symplectic group}. The surface $X$
admits $G$ as group of symplectic automorphisms if and only if
$NS(X)=\mathcal{L}_G$.\end{prop} In general it is an open problem
to understand if for a particular lattice $\mathcal{L}_G$ there
exists a K3 surface $X$ such that $NS(X)=\mathcal{L}_G$.\\
For the lattices $\mathcal{L}^{2d}_{\Z/pZ}$  with $p=2,3,5$ and
the lattices
 $\mathcal{L'}^{2d}_{\Z/pZ,r}$ with $(p,r)=(2,2),
(3,3), (5,5), (7,7)$ the existence is proved in \cite{bert Nikulin
involutions} and in \cite{symplectic prime}. Observe that an even
lattice of signature $(1,n)$, $(n<20)$ is the Néron Severi group
of a K3 surface if and only if it admits a primitive embedding in
$\Lambda_{K3}$ (cf. \cite[Corollary 1.9]{morrison}). The lattices
$\mathcal{L'}^{2d}_{G,r}$ for $G=\Z/4\Z$, $\Z/6\Z$ and
$(\Z/2\Z)^2$ and any $r$ satisfy the conditions of \cite[Theorem
1.14.4]{Nikulin bilinear}, so they admit a primitive embedding in
the K3 lattice. The lattices $\mathcal{L}^{2d}_{\Z/7Z}$, $d\equiv
0 \mod 7$, $\mathcal{L}^{2d}_{G}$, $G=\Z/8\Z$, $(\Z/2\Z)^m,$
$m=3,4$, $(\Z/4\Z)^2$, $\Z/2\Z\times \Z/4\Z$,
$\mathcal{L}^{2d}_{G}$ for $G=(\Z/3\Z)^2$, $\Z/2\Z\times \Z/6\Z$
for $d\equiv 0 \mod 3$ do not admit primitive embeddings in the K3
lattice, indeed a careful analysis shows that the discriminant
group admits a number of generators which would be bigger than the
rank of the transcendental lattice (if a K3 surface would exists)
so this is not possible. The other cases needs a detailed
analysis. We discuss the case $G=\Z/7\Z$ and show the following
proposition (some other cases are described in \cite{alitesi}).
The following lemma is completely trivial, but it is useful for
the next proposition.
\begin{lemma}\label{lemma: sublattices lambdaK3} Let $S$ be a lattice primitively embedded in $\Lambda_{K3}$ and $s\in S$.
Then either $T=\Z s\oplus S^{\perp_{\Lambda_{K3}}}$ is primitively
embedded in $\Lambda_{K3}$, or there exists an overlattice $T'$ of
$T$ such that $T'$ is primitively embedded in $\Lambda_{K3}$ and
$T\hookrightarrow T'$  with finite index. Viceversa if $\Gamma=\Z
L\oplus S^{\perp_{\Lambda_{K3}}}$ or an overlattice with finite
index of $\Gamma$ is primitively embedded in $\Lambda_{K3}$, then
$L\in S$ and so there exists an element of square $L^2$ in $S$.
\end{lemma}

\begin{prop}\label{zeta7}
There exists a K3 surface with Néron Severi group isomorphic to
$\mathcal{L}^{2d}_{\Z/7\Z}$ if and only if $d\equiv 1,2,4 \mod 7$
and a K3 surface with N\'eron--Severi group isomorphic to
$\mathcal{L}^{'2d}_{\Z/7\Z,7}$ if and only if $d\equiv 0\mod 7$.
\end{prop}
\bprf We apply Lemma \ref{lemma: sublattices lambdaK3} in our
case. First we consider the lattice
$\Omega^{\perp}_{\Z/7\Z}$$\simeq$$\left(\begin{array}{cc}4&1\\1&2\end{array}\right)\oplus
U(7)$. Let $d$ be a positive integer, then if
\begin{eqnarray}\label{condb} d\equiv 0,1,2,4\mod 7,
\end{eqnarray}
there exists an element in $\Omega^{\perp}_{\Z/7\Z}$ with square
$2d$. These elements are: $$\begin{array}{ll} (0,0,1,k), \mbox{ if
}d=7k,& (1,1,1,k), \mbox{ if }d=7k+4,\\ (0,1,1,k), \mbox{ if
}d=7k+1,&(1,0,1,k), \mbox{ if }d=7k+2.\end{array}$$ If $d$ does not
satisfy \eqref{condb}, there is no element $L$ of square $2d$.
Indeed the bilinear form on the vector $(p,q,r,s)$ is
$4p^2+2pq+2q^2+14rs$. If $L=(p,q,r,s)$, then $d=2p^2+pq+q^2+7rs$.
Considering this equation modulo 7 we obtain $d=(4p+q)^2$, but 3,
5 and 6 are not square in $\Z/7\Z$. Hence the lattices
$\mathcal{L}^{2d}_{\Z/7\Z}$ with $d\equiv 3,5,6 \mod 7$ do not
admit a primitive embedding in $\Lambda_{K3}$ by Lemma \ref{lemma:
sublattices lambdaK3}. Since the lattice
$\mathcal{L}^{2d}_{\Z/7\Z}$ with $d\equiv 1,2,4 \mod 7$ does not
admit overlattices (Proposition \ref{prop: possible NS of K3 with
G symplectic}), it admits a primitive embedding in
$\Lambda_{K3}$, again by Lemma \ref{lemma: sublattices lambdaK3}.
We observed before that $\mathcal{L}^{2d}_{\Z/7\Z}$ for $d\equiv
0\mod 7$ is not primitively embedded in $\Lambda_{K3}$. On the
other hand, by Lemma \ref{lemma: sublattices lambdaK3},
$\mathcal{L}^{2d}_{\Z/7\Z}$ or an overlattice of
$\mathcal{L}^{2d}_{\Z/7\Z}$ has to be primitively embedded in
$\Lambda_{K3}$. Hence the lattice $\mathcal{L}^{'2d}_{\Z/7\Z,7}$
is primitively embedded in $\Lambda_{K3}$ for each $d\equiv 0\mod
7$.\eprf In the cases where the existence is proved one can give a
description of the coarse moduli space similar to
\cite[Proposition 2.3]{bert Nikulin involutions}, \cite[Corollary
5.1]{symplectic prime}. In fact under this condition, Propositions
\ref{prop: possible NS of K3 with G symplectic}, \ref{prop: if LG
is Ns then G acts symplectically} imply that the coarse moduli
space of algebraic K3 surfaces admitting a certain group of
symplectic automorphisms is the coarse moduli space of the
$\mathcal{L}$-polarized K3 surfaces, for a certain lattice
$\mathcal{L}$ (for a precise definition of
$\mathcal{L}$-polarized K3 surfaces see \cite{dolgachev}).
\begin{rem}{\rm The dimension of the coarse moduli space of the algebraic K3 surfaces admitting $G$ as
group of symplectic automorphisms is $19-\rank(\Omega_G)$.\\
\binf for each group $G$ there exist K3 surfaces such that their
N\'eron--Severi group is either $\mathcal{L}_G^{2d}$ or
$\mathcal{L}_{G,r}^{'2d}$ (cf. \cite{alitesi} and in fact is a
similar computation as in Proposition \ref{zeta7}), so the generic
K3 surface admitting $G$ as group of symplectic automorphisms has
Picard number $1+\rank(\Omega_G)$ and so the moduli space has
dimension $19-\rank(\Omega_G)$. \erem}
\end{rem}
{\it {\bf Acknowledgments.} We warmly thank Bert van Geemen for discussions and for his
constant support during the preparation of this paper. We thank
also Gabriele Nebe for comments and suggestions.}

\newpage

\section{Appendix}\label{appendix}
In Tables 1, 2 we number the fibers according to the
column``singular fibers" (the fiber 1 is the first fiber mentioned
in the column ``singular fibers") and we write the intersections
of the section and the components of the fibers which are not
zero. We denote by $N$ the Nikulin lattice (see \cite{morrison} for
a description).
%For simplicity in table \ref{table3} we denote by $\Omega_r$ the sublattice
%$(H^2(X,\Z)^{\Z/r\Z})^{\perp}$ and we use a similar notation in the case of the groups $(\Z/2\Z)^t$ or $\Z/r\Z \times \Z/s\Z$.

\begin{center}
\section*{Table 1}\label{table1}
\end{center}

{\tiny
%\begin{figure}
\begin{center}
\begin{rotate}{-90}
$
\begin{array}{l|l|c|c|l|c}
G&\mbox{equation}&\mbox{trivial lattice}&\mbox{singular fibers}&\mbox{intersection torsion section}&T_X\\
\hline

\Z/2\Z&\begin{array}{l}y^2=x(x^2+a(\tau)x+b(\tau)),\\\deg(a)=4,\
\deg(b)=8\end{array}& U\oplus A_1^{\oplus
8}&8I_2+8I_1&\begin{array}{l}t_1\mbox{ has order 2 }\\
t_1\cdot C_1^{(i)}=1,\ i=1,\ldots, 8\end{array}&N\oplus U(2)^{\oplus 2}
\\
\hline
\Z/3\Z&\begin{array}{l}y^2=x^3+A(\tau)x+B(\tau)),\\A=\frac{\textstyle
6dc+d^4}{\textstyle 3},~ B=\frac{\textstyle 27c^2-d^6}{\textstyle
3^3}\end{array}& U\oplus A_2^{\oplus
6}&6I_3+6I_1&\begin{array}{l}t_1\mbox{ has order 3 }\\
t_1\cdot C_1^{(i)}=1,\ i=1,\ldots, 6.\end{array}&U\oplus
U(3)\oplus
A_2(-1)^{\oplus 2}\\
\hline

\Z/4\Z&\begin{array}{l}y^2=x(x^2+(e^2(\tau)-2f(\tau))x+f^2(\tau)),\\\deg(f)=4,\
\ \deg(e)=2\end{array}& U\oplus A_3^{\oplus 4}\oplus A_1^{\oplus
2}&4I_4+2I_2+4I_1&\begin{array}{l}t_1\mbox{ has
order 4 }\\
t_1\cdot C_1^{(i)}=1,\ i=1,\ldots, 6.\end{array}&\\
\hline

\Z/5\Z&\begin{array}{l}y^2=x^3+A(\tau)x+B(\tau),\\
A=\frac{\textstyle (-g^4+g^2h^2-h^4-3hg^3+3h^3g)}{\textstyle 3},\\
B=\frac{\textstyle
(g^2+h^2)(19g^4-34g^2h^2+19h^4)}{\textstyle 108}\\
\frac{\textstyle +18hg^3-18h^3g}{\textstyle 108}\\
\end{array}& U\oplus A_4^{\oplus
4}&4I_5+4I_1&\begin{array}{l}t_1\mbox{ has order 5 }\\t_1\cdot
C_1^{(i)}=1,\ i=1,2,\\ t_1\cdot C_2^{(j)}=1,\ j=3,4.
\end{array}&U\oplus U(5)\\
\hline

\Z/6\Z&\begin{array}{l}y^2=x(x^2+(-3k^2(\tau)+l^2(\tau))x+k^3(\tau)(3k(\tau)+2l(\tau))),\\
\deg(k)=\deg(l)=2\end{array}& U\oplus A_5^{\oplus 2}\oplus
A_2^{\oplus 2}\oplus A_1^{\oplus
2}&2I_6+2I_3+2I_2+2I_1&\begin{array}{l}t_1\mbox{ has order 6
}\\t_1\cdot C_1^{(i)}=1,\ i=1,\ldots, 6.\end{array}&U\oplus U(6)
\\
\hline

\Z/7\Z&y^2+(1+\tau-\tau^2)xy+(\tau^2-\tau^3)y=x^3+(\tau^2-\tau^3)x^2
& U\oplus A_6^{\oplus 3}&3I_7+3I_1&\begin{array}{l}t_1\mbox{ has
order 7}\\t_1\cdot C_1^{(1)}=t_1\cdot C_2^{(2)}=1,\\ t_1\cdot
C_3^{(3)}=1.\end{array}&\left[\begin{array}{cc}4&1\\
        1&2\end{array}
        \right]\\
\hline

\Z/8\Z&\begin{array}{l}y^2=x\left(x^2+\left(-2m(\tau)^2n(\tau)^2+\frac{(m(\tau)-n(\tau))^4}{4}\right)x\right)\\
+x\left(m^4(\tau)n^4(\tau)\right),\\
 \deg(m)=1,\ \ \deg(n)=1. \end{array}& U\oplus A_7^{\oplus 2}\oplus
 A_3\oplus A_1
&2I_8+I_4+I_2+2I_1&
\begin{array}{l}t_1\mbox{ has order 8}\\t_1\cdot
C_1^{(i)}=1,\ i=1,3,4,\\ t_1\cdot
C_3^{(2)}=1\end{array}&\left[\begin{array}{cc}2&0\\0&4\end{array}\right]\\
\end{array}
$
%\caption{Table 1}
\end{rotate}
\end{center}
%\end{figure}
}

\newpage

\begin{center}
\section*{Table 2}\label{table2}
\end{center}

{\tiny
\begin{center}
\begin{rotate}{-90}
$
\begin{array}{l|l|c|c|l|c}
G&\mbox{equation}&\mbox{trivial lattice}&\mbox{singular fibers}&\mbox{intersection torsion section}&T_X\\
\hline

\Z/2\Z\times\Z/2\Z&\begin{array}{l}y^2=x(x-p(\tau))(x-q(\tau)),\\\deg(p(\tau))=4,\
\deg(q(\tau))=4\end{array}& U\oplus A_1^{\oplus
12}&12I_2&\begin{array}{l}t_1,u_1\mbox{ have order 2 },\\t_1\cdot
C_1^{(i)}=1,\ i=1,\ldots,8,\\ u_1\cdot
C_1^{(j)}=1,\ j=5,\ldots,12.\end{array}&U(2)^{\oplus 2}\oplus A_1^{\oplus 4}\\
\hline

\Z/2\Z\times\Z/4\Z&\begin{array}{l}y^2=x(x-r^2(\tau))(x-s^2(\tau)),\\\deg(r)=2,\
\deg(s)=2\end{array}& U\oplus A_1^{\oplus 4}\oplus A_3^{\oplus
4}&4I_4+4I_2&\begin{array}{l}
t_1\mbox{ has order 4 },u_1\mbox{ has order 2 } \\
t_1\cdot C_1^{(i)}=1,\ i=1,\ldots,6,\\ t_1\cdot C_0^{(i)}=1, \ i=7,8,\\
u_1\cdot C_2^{(i)}=1,\ i=1,2,\\ u_1\cdot C_0^{(i)}=1, \ i=3,4,\\
u_1\cdot C_1^{(i)}=1,\ i=5,6,7,8.
\end{array}&U(2)\oplus U(4)\\
\hline

\Z/2\Z\times\Z/6\Z&\begin{array}{l}y^2=x[x-(3w(\tau)-z(\tau))(w(\tau)+z(\tau)^3][x-(3w(\tau)\\
+z(\tau))(w(\tau)-z(\tau))^3],\\
\\\deg(w)=1,\
\deg(z)=1\end{array}& U\oplus A_5^{\oplus 3}\oplus A_1^{\oplus
3}&3I_2+3I_6&\begin{array}{l}
t_1\mbox{ has order 6 },u_1\mbox{ has order 2 } \\
t_1\cdot C_1^{(i)}=1,\ i=1,2,4,5,\\
t_1\cdot C_2^{(3)}=t_1\cdot C_0^{(6)}=1,\\
u_1\cdot C_3^{(i)}=1,\ i=1,3,\\
u_1\cdot C_0^{(i)}=1, \ i=2,4,\\
 u_1\cdot C_1^{(i)}=1, \ i=5,6.\end{array}&\left[\begin{array}{cc}6&0\\0&2\end{array}\right]\\
\hline

\Z/3\Z\times\Z/3\Z&\begin{array}{l}y^2=x^3+12x[(\tau^2+1)(\alpha \tau^2+\beta)^3-(\tau^2+1)^4]+\\
2[(\alpha \tau^2+\beta)^6-20(\alpha
\tau^2+\beta)^3(\tau^2+1)^3-8(\tau^2+1)^6],
\\\alpha,\beta\in\C \end{array}& U\oplus A_2^{\oplus
8}&8I_3&\begin{array}{l}t_1, u_1\mbox{ have order 3 }\\t_1\cdot
C_1^{(i)}=1, \ i=1,\ldots,6,\\
t_1\cdot C_0^{(i)}=1, \ i=7,8,\\
u_1\cdot C_0^{(i)}=1, \ i=1,2,\\
u_1\cdot C_1^{(i)}=1, \ i=3,4,7,8,\\
u_1\cdot C_2^{(i)}=1, \ i=5,6. \end{array}&\left[\begin{array}{cc}2&3\\3&0\end{array}\right]\\
\hline

\Z/4\Z\times\Z/4\Z&\begin{array}{l}y^2=x(x-u^2(\tau)v^2(\tau))(x-(1/4)(u^2(\tau)+v^2(\tau))^2),
\\\deg(u)=1,\ \deg(w)=1 \end{array}& U\oplus A_3^{\oplus
6}&6I_4&\begin{array}{l}t_1,u_1\mbox{ have order 4 }\\
t_1\cdot C_1^{(i)}=1,\ i=1,2,3,4,\\ t_1\cdot C_2^{(5)}=t_1\cdot C_0^{(6)}=1,\\
u_1\cdot C_1^{(i)}=1,\ i=3,4,5,6,\\ u_1\cdot C_2^{(1)}=t_1\cdot C_0^{(2)}=1.\\
\end{array}&\left[\begin{array}{cc}4&0\\0&4\end{array}\right]\\

\end{array}
$

\end{rotate}
\end{center}
%\end{figure}
}
\newpage

In Table 3 we give a set of generators for the lattice $\Omega_G$
for all the group $G$ in the list \eqref{formula: symplectic
group} except $G=(\Z/4\Z$), $i=1,2$,  $G=(\Z/2\Z)^i$, $i=2,3,4$,
$G=\Z/2\Z\times \Z/4\Z$, which are presented explicitly in
Propositions \ref{omega and invariant 44}, \ref{omega and
invariant 24 22 4}, \ref{omega and invariant 222,2222}. We always
refer to the elliptic fibration given in Table 2. In each case the
sublattice of the N\'eron--Severi group which is fixed by $G$ is
generated, at least over $\Q$, by the class of the fiber and by
the sum of the classes of the sections.
\begin{center}
\section*{Table 3}\label{table3}
\end{center}

{\tiny
\begin{center}
\begin{rotate}{-90}
$
\begin{array}{c|l}
G&\mbox{generators of the lattice }\Omega_G=(NS(X_G)^G)^{\perp}\\
\hline
\Z/2\Z&b_i=C_1^{(i)}-C_1^{(i+1)},\ i=1,\ldots, 6,\ b_7=F-C_1^{(1)}-C_1^{(2)},\ b_8=s-t.\\
\hline \Z/3\Z& \begin{array}{llll} g_1=s-t_1,&  g_2=s-t_2,&
g_{i+2}=C_1^{(i)}-C_1^{(i+1)},&
  i=1,2,3,4,\\  g_{j+6}=C_1^{(j)}-C_2^{(j+1)},&   j=1,2,3,4,&
 g_{11}=C_2^{(1)}-C_2^{(2)},&  g_{12}=C_0^{(1)}-C_2^{(2)}.\end{array}\\
 \hline
\Z/5\Z&\begin{array}{l}b_1=s-t_1,\ b_2=t_1-t_2,\ b_3=t_2-t_3,\
b_4=t_3-t_4,\ b_5=C_0^{(3)}-C_1^{(2)},\
b_{i+5}=C_1^{(i)}-C_3^{(i)},\ i=1,2,3,\\
b_{j+8}=C_2^{(j)}-C_4^{(j)},\ j=1,2,3,\
b_{h+11}=C_1^{(h)}-C_4^{(h)},\ h=1,2,3,\
b_{k+14}=C_1^{(k)}C_1^{(k+1)},\ k=1,2.\end{array}
\\
\hline \Z/6\Z&\begin{array}{l} b_i=C^{(2)}_i-C^{(2)}_{i+2},\
i=1,\ldots,4,\ b_5=C^{(2)}_2-C^{(1)}_1,\ b_6=C^{(2)}_1-C^{(1)}_1,\
b_7=C^{(1)}_1-C^{(1)}_3,\ b_8=C^{(1)}_1-C^{(1)}_4,\
b_9=C^{(1)}_2-C^{(1)}_4,\ b_{10}=C^{(3)}_1-C^{(4)}_1,\\
b_{11}=C^{(3)}_1-C^{(4)}_2,\ b_{12}=C^{(3)}_2-C^{(4)}_1,\
b_{13}=C^{(5)}_1-C^{(6)}_1,\
b_{14}=C^{(6)}_1-C^{(1)}_1-C^{(1)}_2-C^{(1)}_3,\
b_{15}=C^{(4)}_2-C^{(1)}_1-C^{(1)}_2,\ b_{16}=s-t.
\end{array}\\
\hline \Z/7\Z&\begin{array}{l} b_1=s-t_1,\
b_2=C_0^{(2)}-C_1^{(1)},\ b_{i+2}=t_i-t_{i+1},\ i=1,\ldots, 5,\
b_{j+7}=C^{(1)}_{j}-C^{(1)}_{j+2},\ j=1,\ldots, 4,\\
b_{12}=C_3^{(1)}-C_6^{(1)},\ b_{13}=C_3^{(2)}-C_6^{(2)},\
b_{14}=C_1^{(1)}-C_1^{(2)},\ b_{h+14}=C^{(2)}_{h}-C^{(2)}_{h+2},\
h=1,\ldots, 4.
\end{array}
\\
\hline \Z/8\Z&\begin{array}{llll} b_{1}=t_1-s, &
b_{2}=t_3-t_1,&b_{3}=t_4-t_2,&
b_4=C^{(3)}_3-C^{(1)}_1-C^{(1)}_2,\\
b_{i}=C^{(1)}_{i-2}-C^{(1)}_{i-4},&i=5,6,7,8,&
b_{j}=C^{(1)}_{j-6}-C^{(1)}_{j-4},&j=9,\ldots,13,\\
b_{14}=C_4^{(1)}-C^{(1)}_1,& b_{15}=C_4^{(2)}-C^{(2)}_1,&
b_{16}=C_0^{(2)}-C_1^{(1)},& b_{17}=C_1^{(2)}-C_{1}^{(1)},\
b_{18}=C_1^{(3)}-C_3^{(3)}.
\end{array}
\\
\hline \Z/2\Z\times \Z/6\Z&\begin{array}{llll}
b_1=s-t_1,&b_2=s-u_1,&b_{i+2}=C_1^{(i)}-C_3^{(i)},\
i=1,2,3,&b_{j+5}=C_2^{(j)}-C_4^{(j)},\ j=1,2,3,\\
b_9=C_3^{(1)}-C_5^{(1)},&
b_{10}=C_3^{(3)}-C_5^{(3)},&b_{11}=C_1^{(1)}-C_1^{(2)},&b_{12}=C_1^{(4)}-C_1^{(5)},\\
b_{13}=C_1^{(2)}-C_2^{(3)},&
b_{14}=C_0^{(1)}-C_1^{(2)},&b_{15}=C_2^{(1)}-C_2^{(2)},&b_{16}=C_1^{(1)}-C_2^{(2)},\\
&b_{17}=C_1^{(1)}-C_3^{(3)},&b_{18}=C_1^{(5)}-C_1^{(2)}-C_2^{(2)}-C_3^{(2)}.
\end{array}
\\
\hline \Z/3\Z\times \Z/3\Z&\begin{array}{llll}
b_1=s-t_1,&b_2=s-u_1,&b_{i+2}=C_1^{(i)}-C_1^{(i+1)},&i=1,\ldots 7,\\
b_{j+8}=C_2^{(j)}-C_2^{(j+1)},&j=2,\ldots,
6,&b_{15}=C_1^{(1)}-C_2^{(2)},&b_{16}=C_0^{(1)}-C_2^{(2)}.
\end{array}
\\
\end{array}
$
\end{rotate}
\end{center}
}
\newpage

\addcontentsline{toc}{section}{ \hspace{0.5ex} References}

\end{document}